\documentclass[a4paper,11pt]{article}
\usepackage[english]{babel}
\usepackage[tbtags]{amsmath}
\usepackage{amssymb}
\usepackage{amsfonts}
\usepackage{amsthm}
\usepackage{calrsfs}
\usepackage{array}
\usepackage{bbm}
\usepackage{stmaryrd}
\usepackage{upgreek}
\usepackage[svgnames]{xcolor}
\usepackage{hyperref}
\usepackage{mathtools}

\def\div{\, \mbox{div}\,}
\def\d{\displaystyle}
\def\w{\rho }
\def\ps{\phi(s)}

\numberwithin{equation}{section}
\def\ibint{\int_{-1}^1}
\def\iint{\int_{-1}^1}
\def\grad{\nabla}

\def\nn{\nu_0}

\def\v{{\mathrm{d}}v}

\def\Box{\hfill\rule{2.5mm}{2.5mm}}
\def\t{{\mathrm{d}}\tau}

\def\y{{\mathrm{d}}y}

\def\no{\nonumber}

\def\R{{\mathbb {R}}}

\newcommand\taa{s\ge     -\ln T_0}

\newcommand\RR{{\mathcal R}}
\renewcommand\SS{{\mathcal S}}

\theoremstyle{plain}
\newtheorem{thm}{Theorem}
\newtheorem*{thm*}{Theorem}

\newtheorem{prop}{Proposition}[section]
\newtheorem{cors}{Corollary}[section]
\newtheorem{lem}[prop]{Lemma}

\theoremstyle{remark}
\newtheorem{nb}{Remark}[section]

\usepackage{xcolor}
\usepackage{ulem}

\makeatletter
\def\blfootnote{\xdef\@thefnmark{}\@footnotetext}
\makeatother
\title{\bf 
Blow-up rate for 1-d semilinear wave equations without scaling invariance at characteristic points
}
\date{}

\author{Mohamed  Ali  Hamza\\
{\it \small 
Imam Abdulrahman Bin Faisal University,
 Dammam, 34212, Saudi Arabia}\\
Hatem Zaag\\
{\it \small Universit\'e Sorbonne Paris Nord},\\
{\it \small LAGA, CNRS (UMR 7539), F-93430, Villetaneuse, France}
}

\begin{document}
\maketitle

\begin{abstract}
We investigate the blow-up behavior of solutions to the semilinear wave equation
$$ \partial_t^2 u - \partial_x^2 u = |u|^{p-1}u \ln^a(u^2+2), \quad (x,t)\in \mathbb{R} \times [0,T), $$
with $p > 1$ and $a<0$. While the blow-up rate at non-characteristic points was determined for all $a \in \mathbb{R}$ in previous work, the analysis at characteristic points has remained an open problem. The main obstruction is that known constructions of a Lyapunov functional in similarity variables require several steps and rely on a covering argument, which is not available at characteristic points. Focusing on the regime $a < 0$, we construct a Lyapunov functional directly, without any covering argument. This yields uniform bounds on the solution in similarity variables, and hence the blow-up rate at every point, including characteristic ones. To our knowledge, this is the first work to determine the blow-up rate at characteristic points for an equation  without scaling invariance.
\end{abstract}


\section{Introduction}
This paper is devoted to the study of blow-up solutions for the
following semilinear  wave equation:
\begin{equation}\label{gen}
\left\{
\begin{array}{l}
\partial_t^2 u =\partial_x^2u +|u|^{p-1}u \ln^a (u^2+2),\qquad  (x,t)\in \R\times [0,T),\\
\\
u(x,0)=u_0(x)\in  H^{1}_{loc,u}(\R),\qquad 
\partial_tu(x,0)=u_1(x)\in  L^{2}_{loc,u}(\R),
\end{array}
\right.
\end{equation}
where $u(t):x\in{\R} \rightarrow u(x,t)\in{\R}$, $ p>1$ and $a\in\R$ (our main result concerns the case $a<0$).
The spaces  $L^{2}_{loc,u}(\R)$ and $H^{1}_{loc,u}(\R)$ are  defined by
\begin{equation*}
L^{2}_{loc,u}(\R)=\{u:\R\rightarrow \R/ \sup_{d\in \R}\Big(\int_{|x-d|\le 1}|u(x)|^2dx\Big)<+\infty \},
\end{equation*}
and
\begin{equation*}
H^{1}_{loc,u}(\R)=\{u\in L^{2}_{loc,u}(\R),\partial_x u\in L^{2}_{loc,u}(\R) \}.
\end{equation*}

\medskip

  The Cauchy problem for the  equation ({\ref{gen}}) is well-posed in
$H^{1}_{loc}(\R)\times L^{2}_{loc}(\R)$, locally in
time. This follows from the finite speed
of propagation and the well-posedness in $H^{1} (\R)\times L^{2}(\R)$.
Combining the existence of blow-up solutions of the associated ordinary differential equation with the finite speed of
propagation, one obtains a blow-up solution $u(t)$ of \eqref{gen} which depends non-trivially on the space variable.

\bigskip

In the present paper, we determine the blow-up rate of solutions at a characteristic point for this non-scaling-invariant equation. More precisely, we show
 that all the energy terms in
the similarity variables remain bounded (Theorem~\ref{t1} below). Related questions at non-characteristic points have been addressed in various contexts, ranging from perturbative problems \cite{HZjhde12, HZnonl12, HZkg12, omar1, omar2, H1} to the analysis of equation \eqref{gen} in \cite{HZjmaa2020, HZ2022}. Recently, \cite{RZ} extended these results to nonlinearities involving double-logarithmic perturbations of the form $|u|^{p-1}u \ln^a (\ln(10 + u^2))$ for $a \in \R$. In all those studies, the construction of a Lyapunov functional in 
the similarity variables requires a multi-step procedure and relies on a covering technique, which is inherently restricted to the non-characteristic setting. In this paper we overcome this difficulty when $a<0$ by constructing a Lyapunov functional directly, without any covering argument. 
To the best of our knowledge, the present work is the first to overcome this difficulty for $a<0$, which is essential for a complete understanding of the blow-up dynamics at characteristic points.

\bigskip

If $u$ is an arbitrary blow-up solution of   \eqref{gen}, we define
a 1-Lipschitz curve $\Gamma=\{(x,T(x))\}$
such that
the domain of definition of $u$ is written as
\begin{equation}\label{defdu}
D=\{(x,t)\;|\; t< T(x)\}.
\end{equation}
The justification of this fact can be found in Alinhac
\cite{Apndeta95}, where such a set is referred to as the \textit{maximal influence domain}. 
 In fact, for any $x_0\in \R$, we derive a solution
which is defined in $K_{x_0,t_0}$, the backward light cone
defined by
\[
K_{x_0,t_0}=\{(x,t)\in \R\times[0,t_0)\;|\;|x-x_0|<t_0-t\}.
\]

\medskip

The time $\bar T=\inf_{x\in {\R}}T(x)$ and the curve  $\Gamma$ are called the blow-up time and the blow-up graph of $u$.
A point $x_0$ is a non-characteristic point
if  there exist
\begin{equation}\label{nonchar}
\delta_0\in(0,1)\mbox{ and }t_0<T(x_0)\mbox{ such that }
u\;\;\mbox{is defined on }{\mathcal C}_{x_0, T(x_0), \delta_0}\cap \{t\ge t_0\}
\end{equation}
where ${\cal C}_{\bar x, \bar t, \bar \delta}=\{(x,t)\;|\; t< \bar t-\bar \delta|x-\bar x|\}$. We denote by $\RR$ (resp. $\SS$) the set of non-characteristic (resp. characteristic) points.

\bigskip
For the standard semilinear wave equation with a pure power nonlinearity,
\begin{equation}\label{NLW}
\partial^2_{t} u = \Delta u + |u|^{p-1}u, \quad (x,t) \in \mathbb{R}^n \times [0,T),
\end{equation}
the blow-up rate of blow-up solutions was  studied by Merle and Zaag in \cite{MZajm03, MZimrn05, MZma05}, for $p > 1$ when $n = 1$, and in the subconformal and conformal regimes $1 < p \le \frac{n+3}{n-1}$ when $n \ge 2$.
 Their approach relies heavily on the analysis of the equation in  similarity variables. More precisely, for any given $x_0 \in \mathbb{R}^n$ and $T_0 \in (0, T(x_0)]$, we introduce the similarity variables defined by:
$$ y = \frac{x-x_0}{
T_0-t} ,\ \  s = -\ln(T_0- t),\quad 
\textrm{and} \qquad w_{x_0,T_0} (y, s)= (T_0-t)^{\frac2{p-1}}u(x, t).$$
 From \eqref{NLW}, the function $w_{x_0,T_0}$ (hereafter denoted by $w$ for simplicity) satisfies the following equation for all $y \in B$ and $s \ge -\ln T_0$, where $B \equiv B(0,1)$ denotes the unit ball in $\mathbb{R}^n$:
\begin{equation*}\label{w0}
\partial_{s}^2w=\frac{1}{\rho}\div(\rho \grad w-\rho (y \cdot \nabla w)y)
-\frac{2p+2}{(p-1)^2}w
-\frac{p+3}{p-1}\partial_s w
-2y \cdot \nabla \partial_{s}w+
|w|^{p-1}w,
\end{equation*}
where $\rho (y)=(1-|y|^2)^{\alpha}$, with
$\alpha =\frac{2}{p-1}-\frac{n-1}2$.
The analysis is fundamentally based on the construction of a suitable Lyapunov functional in similarity variables and a precise treatment of its dissipative properties. This approach is supplemented by sharp Sobolev embeddings and critical Gagliardo-Nirenberg inequalities to control the nonlinear term. In particular, one derives the following uniform bound on time averages:
\begin{equation}
  \int_s^{s+1}\int_B \Big((\partial_{s}w_{x_0})^2+|\grad w_{x_0}|^2-(y \cdot \nabla w_{x_0})^2+|w_{x_0}|^{p+1} \Big)\rho \y\t\le K_0,
\end{equation}
where,  $w_{x_0}$ denotes $w_{x_0, T(x_0)}$,
and
where 
$K_0 = K_0(n, p, T(x_0), \|(u_0, u_1)\|_{H^1 \times L^2})$  is a positive constant.

\medskip

If $x_0$ is a non-characteristic point, one can employ a geometric covering technique to establish that the solution in similarity variables remains bounded. Specifically, there exist constants $\varepsilon_0 = \varepsilon_0(n, p) > 0$ and $\tilde K_0 = \tilde K_0(n, p, T(x_0), \|(u_0, u_1)\|_{H^1 \times L^2})$ such that
\begin{equation*}
0<\varepsilon_0\le \|w_{x_0}(s)\|_{H^{1}(B)}+ \|\partial_s
w_{x_0}(s)\|_{L^{2}(B)} \le \tilde K_0, \qquad \forall s\ge -\ln T(x_0).
\end{equation*}
More precisely, in the  $u(x,t)$ setting, 
if $x_0$ is a non-characteristic point,    and $ p \in (1, \frac{n+3}{n-1}]$, the blow-up rate is governed by the dynamics of the associated ordinary differential equation (ODE). This behavior is categorized as a Type I blow-up (or the ODE regime); for equation \eqref{gen}, the corresponding rate is $\psi_T$ defined in \eqref{psi} below. Conversely, any blow-up whose rate deviates from this scale is called Type~II.

\medskip

In contrast, the analysis at characteristic points requires a more refined treatment of the energy functional. Regarding the one-dimensional case, Merle and Zaag \cite{MZjfa07} established 
the following sharp uniform estimate, valid at every point, including characteristic ones:
\begin{equation}\label{jfa0}
 \iint \Big((\partial_{s}w)^2+(1-y^2)(\partial_y w)^2+|w|^{p+1} \Big)(1-y^2)^{\frac2{p-1}} \y\le \hat K_0,  \qquad \forall s\ge -\ln T(x_0).
\end{equation}

In view of the estimate \eqref{jfa0} established in \cite{MZjfa07} for the one-dimensional setting, further refinements were obtained in \cite{MZjfa07, MZcmp08, MZajm11, MZisol10, CZcpam13}, ending by proving the resolution in solitons at blow-up. These results were adapted to the radial case away from the origin in \cite{MZbsm11} and further extended to the Klein-Gordon equation and other damped lower-order perturbations in \cite{HZkg12}. In later works \cite{MZods, MZods14}, Merle and Zaag extended their analysis to higher dimensions under subconformal conditions, proving that the explicit self-similar solution is stable with respect to perturbations of the blow-up point and initial data.

\bigskip

We return to equation \eqref{gen}. Before our contributions \cite{HZjmaa2020, HZ2022}, it was  natural to  expect that the blow-up rate for the nonlinear wave equation \eqref{gen} in the non-characteristic case is governed by the dynamics of the associated ordinary differential equation:
\begin{equation}\label{v}
v_T'' (t)=  |v_T (t)|^{p-1}v_T(t) \ln^{a}(v_T^2(t)  +2), \quad v(T)=\infty.
\end{equation}
Note that    the solution of \eqref{v} satisfies 
 \begin{equation}\label{equivv0}
 v_T(t) \sim \kappa_{a} \psi_T(t), \text{ as } t \to\ T,\quad \textrm{ where}\quad  \kappa_{a} =  \left(\frac{2^{1-2a}(p+1)}{(p-1)^{2-a}} \right)^{\frac1{p-1}},
 \end{equation}
and
 \begin{equation}\label{psi}
 \psi_T(t)=(T - t)^{-\frac{2}{p-1}} (-\ln (T - t))^{-\frac{a}{p-1}},
 \end{equation}
(see Lemma A.2  in \cite{HZjmaa2020} for the  proof of \eqref{equivv0}).

 \medskip

As a matter of fact, this blow-up rate was first rigorously established for equation \eqref{gen} in the non-characteristic case in \cite{HZjmaa2020, HZ2022}. Indeed, we demonstrated that the blow-up is of Type I for any $a \in \mathbb{R}$, provided that the exponent is subconformal. Very recently, these results were extended to the conformal regime, i.e., $p = \frac{n+3}{n-1}$ for $n \ge 2$ and $a < 0$, in \cite{Hamza2026}.
The same persistence of the Type~I rate was shown in \cite{RZ} for the double-logarithmic nonlinearity $|u|^{p-1}u\ln^a(\ln(10+u^2))$, $a\in\R$, in the subconformal regime.

 \medskip

While previous investigations into blow-up rates for equation \eqref{gen} have been confined to non-characteristic points, the present work extends these results to include characteristic points.
We avoid the covering argument, which is not available at characteristic points, and instead construct
a Lyapunov functional by employing refined estimates for the time derivative of the energy.
Since these estimates rely on the sign condition $a < 0$, this assumption is required in our main results.

Throughout the paper, unless otherwise specified, $u$
denotes a solution of \eqref{gen} with blow-up graph
$\Gamma = \{ (x, T(x)) : x \in \mathbb{R} \}$, and $x_0\in \R$
is fixed.  Since \eqref{gen} is invariant under time
translations, we may shift the initial time to $t_0(x_0)$, so
that the remaining lifespan $T(x_0) - t_0(x_0)$ is   as small as we wish; the
new initial data are the values of $u$ at $t=t_0(x_0)$. After this shift, $T(x_0)$
denotes the remaining lifespan, and we may assume that
\begin{equation}\label{ss00}
T(x_0)\le T_*,
\end{equation}
where $T_*\in(0,1)$ is a small
constant depending only on $p$ and $a$, fixed once and for all so that every inequality
below requiring $s$ to be large holds for all $s\ge-\ln T_*$. In particular
$-\ln T(x_0)\ge -\ln T_*\gg1$.

\medskip

Under  condition \eqref{ss00}, we introduce the following adapted similarity variables, defined for any
$x_0 \in \R$ and $T_0$ such that $0<T_0\le T(x_0)$, by
\begin{equation}\label{scaling}
y=\frac{x-x_0}{T_0-t}, \quad s=-\ln(T_0-t), \quad
u(x,t)=\psi_{T_0}(t)\, w_{x_0,T_0}(y,s).
\end{equation}
From (\ref{gen}), the  function $w_{x_0,T_0}$  (we write $w$ for
simplicity) satisfies the following equation for all $y\in (-1,1)$ and
$s\ge  -\ln T_0$,
\begin{align}\label{A}
\partial_{s}^2w&=\frac{1}{\rho}\partial_y(\rho(1-y^2)\partial_y w)+\frac{2a}{(p-1)s}
y\partial_y w-\frac{2p+2}{(p-1)^2}w+\gamma(s)w\nonumber\\
&-\Big(\frac{p+3}{p-1}-\frac{2a}{(p-1)s}\Big)\partial_s w
-2y \partial^2_{ys}w+\frac{1}{s^a}   |w|^{p-1}w\ln^a(\phi^2w^2+2),
\end{align}
where   
$\rho (y)=(1-y^2)^{\alpha}$, $\alpha=\frac{2}{p-1},$  
\begin{equation}
   \gamma(s)=\frac{a(p+3)}{(p-1)^2s}-\frac{a(p+a-1)}{(p-1)^2s^2},\label{defgamma}
\end{equation}
and 
\begin{equation}
\ps =e^{\frac{2s}{p-1}}s^{-\frac{a}{p-1}}.\label{defphi}
\end{equation}

In the new variables $w(y,s)$, deriving the behavior of $u$ as $t\to T_0$ is equivalent to deriving the behavior of $w$ as $s\to+\infty$.  Moreover, if $T_0=T(x_0)$, then we simply write $w_{x_0}$
instead of $w_{x_0,T(x_0)}$.

\medskip

This equation will be studied in the space
\begin{equation}\label{H}
\mathcal{H} = \left\{ q=(q_1,q_2) \in H_{loc}^1 \times L_{loc}^2(-1, 1) \mid \|q\|_{\mathcal{H}}^2 \equiv \int_{-1}^1 (q_1^2 + (q_1')^2(1 - y^2) + q_2^2) \rho \, dy < +\infty \right\},
\end{equation}
which is the energy space for $w$. Note that $\mathcal{H} = \mathcal{H}_0 \times L_\rho^2$ where   $L^2_\rho=L^2((-1,1),\rho\,\mathrm{d}y)$ and

\begin{equation}\label{H0}
\mathcal{H}_0 = \left\{ r \in H_{loc}^1(-1, 1) \mid \|r\|_{\mathcal{H}_0}^2 \equiv \int_{-1}^1 (r'^2(1 - y^2) + r^2) \rho \, dy < +\infty \right\}.
\end{equation}

\medskip

The key ingredient is a Lyapunov functional for \eqref{A}, constructed directly in Section~\ref{sec2}. Combined with the Merle--Zaag strategy for the pure power case, it yields our main result.
\begin{thm}\label{t1}
{\bf {(Blow-up rate for equation \eqref{gen})}}. 
Let $p>1$ and  $u$ be a solution to \eqref{gen} with blow-up graph $\Gamma:\{x\mapsto T(x)\}$. Suppose $x_0$ is a characteristic point satisfying the conditions \eqref{ss00} and that  $a<0$. Then, 
 for all
 $s\ge 2-\ln T(x_0)$,
\begin{equation}\label{30}
  \int_{-1}^1 \Big((\partial_{s}w_{x_0})^2+(\partial_y w_{x_0})^2(1-y^2)+w_{x_0}^2+\frac1{s^{a}}|w_{x_0}|^{p+1}\ln^a (2+\phi^2w_{x_0}^2) \Big)\rho(y) \y\le K,
\end{equation}
where $w_{x_0}=w_{x_0,T(x_0)}$ is defined in (\ref{scaling}), 
and 
$K=K(p, a, T(x_0),\|(u(0),\partial_tu(0))\|_{
H^{1}\times
L^{2}(I(x_0,T(x_0))} )$, with $I(x_0,T(x_0))=(x_0-T(x_0),x_0+T(x_0))$.\\
Furthermore,  when $ s \to \infty$
we have
\begin{align}
\int_{-1}^{1} \Big[\tfrac{1}{2}(\partial_y w_{x_0}(y,s))^2 (1 - y^2) + \tfrac{1}{2}(\partial_s w_{x_0}(y,s))^2 + \tfrac{p+1}{(p-1)^2}\, w_{x_0}(y,s)^2 \Big] \rho \, dy &\to \frac{p+1}{p-1}\, E_{\infty}, \label{33}\\
\frac{1}{s^a} \int_{-1}^{1} |w_{x_0}(s)|^{p+1} \ln^a (2+\phi^2 w_{x_0}^2)\, \rho(y)\, dy &\to \frac{2p+2}{p-1}\, E_{\infty}, \label{34}
\end{align}
for some  $E_{\infty}\ge0$.
\end{thm}

\begin{nb}
The nonnegative real number $E_{\infty}$ is the limit, as $s\to\infty$, of the Lyapunov
functional $L$ introduced in \eqref{LL}. Its existence and nonnegativity follow from
Proposition~\ref{lyap}, which shows that $L$ is monotone and bounded below by zero.
For $a=0$, Theorem~\ref{t1} recovers the result of \cite[Proposition~3.5, p.~66]{MZjfa07}.
\end{nb}

\begin{nb}
In Theorem~\ref{t1}, we choose to assume that $x_0 $ is a characteristic point (rather than an arbitrary point, since the result is already known for non-characteristic points, with the primary case of interest being that of a characteristic point).
\end{nb}

\begin{nb}
Theorem~\ref{t1} extends the known blow-up results of \cite{MZjfa07} for the standard semilinear wave equation \eqref{NLW} to the more general equation \eqref{gen}, which is not invariant under scaling. To the best of our knowledge, it provides the first  determination of the blow-up rate at characteristic points for a non-scaling-invariant equation.
\end{nb}

\begin{nb}
As noted above,  the time-translation invariance of \eqref{gen} allows us to shift the initial time to some $t_0(x_0)<T(x_0)$ such that the remaining lifespan $T(x_0)-t_0(x_0)$ is arbitrarily small. Thus, without loss of generality, we may assume \eqref{ss00}, which implies $-\ln T(x_0)\gg 1$.
 In this configuration, the solution values at $t=t_0(x_0)$ serve as the effective initial data for our analysis. Naturally,  
 by the local well-posedness of \eqref{gen}, these data depend continuously on the original initial data at $t=0$.
\end{nb}

\begin{nb}
The same strategy is expected to extend to the loglog nonlinearity $|u|^{p-1}u\ln^a(\ln(10+u^2))$ of Roy and Zaag \cite{RZ} for $a<0$, so their results also hold at characteristic points. For $a>0$, one expects the same conclusion, but this remains an open question.
\end{nb}



\bigskip

Throughout this paper,
$C$  denotes a  generic positive constant
 depending only on $p$  and $a,$  which may vary from line to line.
 In addition, we  will use $K_1,K_2,K_3...$  as    positive constants
 depending only on $p, a,$  and initial data,   which
 may also vary from line to line.  

\bigskip

The remainder of this paper is structured as follows. In Section \ref{sec2}, we show that the functional $L(w(s),s)$ defined in \eqref{LL} serves as a Lyapunov functional for equation \eqref{A}, leading to the derivation of Proposition \ref{lyap}. This result constitutes a pivotal step in determining the blow-up rate for \eqref{A}. 
In Section \ref{section3}, we establish the boundedness of space-time averages and prove Theorem \ref{t1}. 
\section*{Acknowledgment}
Hatem Zaag wishes to thank 
Pierre Rapha\"el
and the ”SWAT”
ERC project for their support.
\section{A positive Lyapunov functional for the solution of (\ref{A}) at a characteristic point}\label{sec2}
This section is devoted to the statement and the proof of Proposition \ref{lyap}. It is divided  into three subsections.
\begin{itemize}
\item Energy estimates: 
 In the first subsection, we derive the foundational energy estimates by testing \eqref{A} with $\partial_s w \rho$ and $w \rho$. A central feature of this analysis is the refined estimate obtained for $a < 0$, which  is
 crucial for the construction of the Lyapunov functional.
\item Lyapunov functional:
The second subsection focuses on the construction of a Lyapunov functional through a suitable combination of the previously derived energy estimates.
\item Positivity: Finally, we establish the positivity of this functional, concluding the proof of the proposition.
\end{itemize}
We start by introducing the following functionals, which are crucial:
\begin{eqnarray}
  E_0(w(s),s)&=&\iint \Big(\frac{1}{2}(\partial_{s}w)^2+\frac{1}{2}(\partial_y w)^2(1-y^2) +\big(\frac{p+1}{(p-1)^2}-\frac12  \gamma (s)\big)w^2\nonumber\\
  &&-e^{-\frac{2(p+1)s}{p-1}}s^{\frac{2a}{p-1}}   f(\phi w)\Big)\rho \y,\label{E0}\\
  E_1(w(s),s)&=&-\frac1{s\sqrt{s}}\iint w\partial_{s}w \rho \y,\label{E1}\\
E(w(s),s)&=&E_0(w(s),s)+ E_1(w(s),s),\label{EE}
\end{eqnarray}
where $\gamma(s)$ and $\phi(s)$  are defined in \eqref{defgamma} and \eqref{defphi}, respectively,
and $f$ is given by 
\begin{equation}\label{defF}
f(u)=\int_0^u|v|^{p-1}v\ln^a(2+v^2)dv.
\end{equation}
Finally, we introduce the following candidate for a Lyapunov functional
\begin{align}
L(w(s),s)=&\exp\left(\frac{p+3}{\sqrt{s}}\right) E(w(s),s)+e^{-s}.\label{LL}
\end{align}
We shall show that $L(w(s),s)$ is non-increasing and nonnegative. More precisely, we have the following:

\begin{prop}\label{lyap}
\textbf{(A positive Lyapunov functional for \eqref{A}.)}
Let $a < 0$ and $x_0 \in \mathbb{R}$. Let $u$ be a solution to \eqref{gen} with blow-up graph
$\Gamma: \{x \mapsto T(x)\}$, and let $T_0 \in (0, T(x_0)]$. 
Assume that \eqref{ss00} holds (here $x_0$ may be characteristic or not), and let
$w = w_{x_0,T_0}$ be the scaling function defined in \eqref{scaling}. Then, for all $s' \ge s \ge -\ln T_0$, the following inequality holds:
\begin{equation}\label{L01}
L(w(s'), s') - L(w(s), s) \le - \frac{2}{p-1} \int_{s}^{s'} \int_{-1}^1 (\partial_{\tau}w)^2 \frac{\rho}{1-y^2} \, \mathrm{d}y \, \mathrm{d}\tau.
\end{equation}
Moreover, for all $s \ge -\ln T_0$, we have
\begin{equation}\label{posi}
L(w(s), s) \ge 0.
\end{equation}
\end{prop}
\begin{nb}
When $x_0$ is a non-characteristic point, this result was already proved for every $a\in\mathbb{R}$ in our previous works \cite{HZjmaa2020, HZ2022}, by means of the standard covering technique. For $a<0$, the proof of Proposition~\ref{lyap} below is direct and does not use any covering argument. Consequently, it holds at every point $x_0\in\mathbb{R}$, whether characteristic or non-characteristic.
\end{nb}

\begin{nb}
The assumption $a<0$ is crucial for the construction of the Lyapunov functional $L(w(s),s)$: it allows us to refine the estimates on the energy $E_0(w(s),s)$, a refinement that is not available when $a>0$.  Consequently, the extension of our results to the case $a > 0$ remains an open question.
\end{nb}

\subsection{Classical energy estimates}\label{sec21}
In this subsection, we present the preliminary lemmas necessary for the construction of a Lyapunov functional for equation \eqref{A}. We begin by establishing an identity for the functional $E_0(w(s), s)$, together with a bound on its remainder, in the following lemma:
\begin{lem}\label{lem21} Let $a<0$. For all  $s'\ge s\ge -\ln T_0$, we have 
 \begin{align}\label{Lem21}
 E_0(w(s'),s')-
E_0(w(s),s)=\!
-\frac{4}{p-1}\int_s^{s'}\int_{-1}^1(\partial_{s} w)^2\frac{\rho}{1-y^2}\y\t +\int_s^{s'}\Sigma_0(\tau)\t,
\end{align}
where $\Sigma_{0}(s)$ satisfies
    \begin{eqnarray}\label{Lem210}
\Sigma_{0}(s)\leq \frac{1}{p-1}\int_{-1}^1(\partial_{s} w)^2\frac{\rho}{1-y^2}\y+
\frac{C_0}{s^{a+\frac74}}\ibint |w|^{p+1}\ln^a(\phi^2w^2+2)\w\y\no\\
+ \frac{C_0}{s^2}\iint (\partial_y w)^2(1-y^2)\rho {\mathrm{d}}y
+ \frac{C_0}{s^2}\iint w^{2} \w {\mathrm{d}}y
+C_0 e^{-2s}.
\end{eqnarray}
\end{lem}
{\it Proof}: 
By multiplying \eqref{A} by $\partial_{s} w$ and integrating over $(-1,1)$, an integration by parts followed by integration over $[s, s']$ yields \eqref{Lem21}, where
\begin{align}\label{estim1}
\Sigma_{0}(s)= &
\underbrace{
\frac{a}{(p+1)s^{a+1}}\ibint |w|^{p+1}\ln^a(\phi^2w^2+2)
\rho\y+\frac{2p+2}{p-1}
e^{-\frac{2(p+1)s}{p-1}}s^{\frac{2a}{p-1}}\ibint f_1(\phi w)
\rho \y}_{\Sigma_0^1(s)}\no\\
&+\underbrace{\frac{2p+2}{p-1}
e^{-\frac{2(p+1)s}{p-1}}s^{\frac{2a}{p-1}}\ibint f_2(\phi w)
\rho \y}_{\Sigma_0^2(s)}\no\\
&\underbrace{-\frac{2a}{(p-1)s}
e^{-\frac{2(p+1)s}{p-1}}s^{\frac{2a}{p-1}}\ibint\big( f_1(\phi w)+f_2(\phi w)\big)
\rho \y}_{\Sigma_0^3(s)}\underbrace{-{\frac{\gamma' (s)}2\ibint w^2\rho\y}}_{\Sigma_0^4(s)}\no\\
&+\underbrace{\frac{2a}{(p-1) s} \int_{-1}^1\partial_{s} w 
 y\partial_y w  \w {\mathrm{d}}y}_{\Sigma_0^5(s)}+\underbrace{\frac{2a}{(p-1) s} \int_{-1}^1(\partial_{s} w )^2
   \w {\mathrm{d}}y}_{\Sigma_0^6(s)},
\end{align}
where  $f_1$, and $f_2$ are given by 
\begin{equation}
f_1(x)= -\frac{ 2a} {(p+1)^2}{| x|^{p+1}}\ln^{{a-1}}(x^2+2 ),\label{defF2}
\end{equation}
and 
\begin{equation}\label{defF123}
f_2(x)=f(x)
-\frac{ 1} {p+1}{| x|^{p+1}}\ln^{a}(x^2+2 )-
f_1(x).
\end{equation}
The term $\Sigma_0^{1}(s)$ is the dominant contributor to the energy estimate, while the remaining terms, $\Sigma_0^{i}(s)$ for $i \in \{2, \dots, 6\}$, are of lower order. We now proceed to derive the necessary bounds for $\Sigma_0^{i}(s)$ for all $i \in \{1, \dots, 6\}$.
By   using the definition of $f_1$, we  show that
\begin{align}\label{re}
\Sigma_0^1(s)= \frac{a}{(p+1)s^{a+1}}\iint  {|  w|^{p+1}}\ln^{{a-1}}(2+\phi^2 w^2  )\Big(\ln (\phi^2w^2+2 )-\frac{4s}{p-1} \Big)\rho \y.
 \end{align}
  For $ \taa$, we introduce the following partition of $B=(-1,1)$
  \begin{equation}\label{partition}
 B_1 (s)=\{y \in B\,\,|\,\,  w^2(y,s)\leq  s^{\frac{2a-4}{p-1}}\}\,\,{\rm and }\,\,  
 B_2(s)=\{y \in B \,\,|\,\,  w^2(y,s)\ge  s^{\frac{2a-4}{p-1}}\}.
\end{equation}
Accordingly, we write  
$\Sigma_0^{1}(s)=\Phi_1(s)+\Phi_2(s)$,  where for $i \in \{1, 2\}$,
\begin{align}
\Phi_i(s)=&
\frac{a}{(p+1)s^{a+1}}\int_{B_i(s)}  {|  w|^{p+1}}\ln^{{a-1}}(2+\phi^2 w^2  )\Big(\ln (\phi^2w^2+2 )-\frac{4s}{p-1} \Big)\rho\y.
\label{i12}
\end{align}
On the one hand, by using the fact $a<0$ (so that $\ln^{a-1}(2+\phi^2w^2)\le (\ln 2)^{a-1}$), the definition of the set $B_1(s)$ which gives $|w|^{p-1}\le s^{a-2}$ on $B_1(s)$,  the bound $|\ln(\phi^2w^2+2)-\frac{4s}{p-1}|\le Cs$  on $B_1(s)$,  the expression of $\phi(s)$ in \eqref{defphi}, and integrating over $B_1(s)$, we get, for all $\taa$,
\begin{equation}\label{ia}
\Phi_1(s) \le \frac{C}{s^2}\iint   w^2  \rho \y.
\end{equation}
By using  the definition of the $\phi$ given by \eqref{defphi}, we write the identity
\begin{equation}\label{i02}
\ln (2+\phi^2 w^2)-\frac{4s}{p-1}
 =\ln (2\phi^{-2}+ w^2)-\frac{2a}{p-1}\ln s.
\end{equation}
Furthermore, utilizing \eqref{i02} along with the definition of the set $B_2(s)$ in \eqref{partition}, we obtain for all
 $\taa$, 
\begin{equation}\label{i03}
\ln (2+\phi^2 w^2)-\frac{4s}{p-1}
 \ge - C\ln s.
\end{equation}
Using again the definition of the set $B_2(s)$ in \eqref{partition} together with the hypothesis \eqref{ss00}, we can write that for all $\taa,$ if $y \in B_2(s)$, then
\begin{equation}\label{i04}
\ln(\phi^2w^2+2)\ge \ln(\phi^2s^{\frac{2a-4}{p-1}})\ge \frac{4s}{p-1}-\frac{4 }{p-1}\ln s\ge \frac{3s}{p-1}.
\end{equation}
Consequently, combining \eqref{i03} and \eqref{i04}, we obtain    for all  $\taa$, 
\begin{align}\label{ib}
\Phi_{2}(s)\le 
\frac{C\ln s}{s^{a+2}}\int_{-1}^1  {|w|^{p+1}}\ln^{{a}}(2+\phi^2 w^2  )\rho \y.
\end{align}
Adding  \eqref{ia} and \eqref{ib}  and  using  the fact 
that
$
\Sigma_0^1(s)=\Phi_1(s)+\Phi_2(s)$,   we get
\begin{equation}\label{j1}
\Sigma_0^1(s)\le 
\frac{C}{s^{a+\frac74}}\iint   {|  w|^{p+1}}\ln^{{a}}(2+\phi^2 w^2 )\rho \y+\frac{C}{s^2}\iint   w^2  \rho \y.
\end{equation}

We next address the remaining terms $\Sigma_0^{i}(s)$ for $i \in \{2, \dots, 6\}$. These terms are of lower order and can be bounded straightforwardly. 
 In what follows, we establish the required estimates for $\Sigma_0^{i}(s)$ for all $i \in \{2, \dots, 6\}$.

\par

First, it follows from \eqref{equiv2} and \eqref{equiv3} that for all $z\in\R$ and $s\ge 1$,
\begin{equation}\label{15dec1}
\frac{1}{s}| f_1(\phi z)|+| f_2(\phi z)|\le   C+\frac{C}{s^2}|\phi z|^{p+1} \ln^a(2+\phi^2 z^2).
\end{equation}
Consequently, \eqref{15dec1} yields for all $\taa$,
\begin{equation}\label{j2}
\Sigma_0^2(s)+\Sigma_0^3(s)\le 
\frac{C}{s^{a+2}}\iint |w|^{p+1}\ln^a(\phi^2w^2+2)\rho \y+  C e^{-2s}.
\end{equation}
Next, recalling the definition of $\gamma(s)$ in \eqref{defgamma}, we have the bound $\vert{}\gamma'(s)\vert{}\le \frac{C}{s^2}$, which implies
\begin{equation}\label{j3}
\Sigma_0^4(s)\leq \frac{C}{s^2}\int_{-1}^1 w^2\rho {\mathrm{d}}y.
\end{equation}
Moreover, the final two terms satisfy
\begin{equation}\label{j4}
\Sigma_{0}^5(s) +\Sigma_{0}^6(s)\le
\frac{1}{p-1}\int_{-1}^1(\partial_{s} w)^2 \frac{\rho}{1-y^2} {\mathrm{d}}y+ \frac{C}{s^2}\int_{-1}^1 (\partial_y w)^2(1-y^2)\rho {\mathrm{d}}y.
\end{equation}
Adding the inequalities  
\eqref{j1},  \eqref{j2}, \eqref{j3},   and \eqref{j4}, yields  \eqref{Lem210}, which ends the proof of Lemma \ref{lem21}.
\Box

 We are now going to prove the following estimate for the functional $E_1(w(s),s)$.
\begin{lem}\label{lem22} 
For all  $s'\ge s\ge -\ln T_0$, we have
\begin{eqnarray}\label{Lem22}
E_1(w(s'),s')-E_1(w(s),s)\!
\!\!
\!&=&\!\!\!-\int_s^{s'}\frac1{\tau^{\frac32}} \int_{-1}^1 (\partial_{s} w)^2\rho {\mathrm{d}}y\t+\int_s^{s'}\frac1{\tau^{\frac32}} 
\int_{-1}^1(1-y^2)(\partial_y w)^2 \rho {\mathrm{d}}y\t \no\\
&&-\int_s^{s'}\!\!\frac{ 1}{\tau^{a+\frac32}}\!\ibint |w|^{p+1}\ln^a(2+\phi^2 w^2)\rho \y\t\no \\
&&\!\! \!\!\!+\int_s^{s'}\!\frac{2p+2}{(p-1)^2\tau^{\frac32}} \int_{-1}^1\!w^2 \rho {\mathrm{d}}y\t+\int_s^{s'}\Sigma_1(\tau)\t,
\end{eqnarray}
where $\Sigma_{1}(s)$ satisfies
\begin{equation}\label{Lem220}
\Sigma_1(s)\le
\frac{C}{s}\int_{-1}^1(\partial_{s} w)^2 \frac{\rho}{1-y^2} {\mathrm{d}}y+
  \frac{C}{s^\frac74}
\int_{-1}^1(\partial_y w)^2(1-y^2)\rho  {\mathrm{d}}y+\frac{C}{s^\frac{7}4}
\int_{-1}^1 w^2\rho  {\mathrm{d}}y.
\end{equation}
\end{lem}
{\it Proof}: Note that $E_1(w(s),s)$ is a differentiable function of $s$. Using equation \eqref{A} and integrating over $[s,s']$, we obtain \eqref{Lem22}, where
\begin{eqnarray}
\label{l21}
 \Sigma_1(s)
&=&\underbrace{- \frac2{s\sqrt{s}}\int_{-1}^1\partial_{s} w y \partial_y w\rho {\mathrm{d}}y}_{ \Sigma^1_1(s)}+\underbrace{ \frac{8}{(p-1)s\sqrt{s}} \int_{-1}^1\partial_{s} w w\frac{y^2\rho}{1-y^2} {\mathrm{d}}y}_{ \Sigma^2_1(s)}
\no\\
&&\underbrace{ (\frac{a}{(p-1)s}
-\gamma (s))\frac1{s\sqrt{s}} \int_{-1}^1w^2 \rho {\mathrm{d}}y
+(\frac{5-p}{p-1}+\frac3{2s}-\frac{2a}{(p-1)s})\frac1{s\sqrt{s}} \int_{-1}^1 w\partial_{s}w \rho \y } _{ \Sigma_1^3(s)}\no\\
&&\underbrace{-\frac{4a}{(p-1)^2s^{\frac52}}\int_{-1}^1w^2 \frac{y^2\rho}{1-y^2} {\mathrm{d}}y} _{ \Sigma_1^4(s)}.
\end{eqnarray}
By applying the basic inequality $2ab \le a^2 + b^2$, we obtain the following estimates:
\begin{equation}\label{bj1}
 \Sigma^1_1(s)\le
\frac{1}{s^{\frac54}}\int_{-1}^1(\partial_{s} w)^2 \frac{\rho}{1-y^2} {\mathrm{d}}y+\frac1{s^{\frac74}}\int_{-1}^1 (\partial_y w)^2(1-y^2)\rho {\mathrm{d}}y,
\end{equation}
and
\begin{equation}\label{bj2}
 \Sigma^2_1(s)\le
\frac{C}{s}\int_{-1}^1(\partial_{s} w)^2 \frac{\rho}{1-y^2} {\mathrm{d}}y+\frac{C}{s^{2}} \int_{-1}^1 w^2\frac{y^2\rho}{1-y^2} {\mathrm{d}}y.
\end{equation}
Furthermore, invoking the estimate $ |\gamma(s)| \le \frac{C}{s}$, we deduce
\begin{align}\label{bj3}
 \Sigma^3_1(s)\le
\frac{C }{s }\int_{-1}^1 (\partial_{s}w)^2\rho  \y+\frac{C }{s^2 }\int_{-1}^1 w^2\rho  \y.
\end{align}
Let us recall from  \eqref{hardy1}, the following Hardy-type inequality
\begin{equation}\label{AHardy}
 \iint w^2\frac{y^2\rho}{1-y^2}{\mathrm{d}}y\leq C\iint (\partial_y w)^2(1-y^2)\rho {\mathrm{d}}y+C\iint w^2\rho {\mathrm{d}}y.
\end{equation}
Thus, it follows from  \eqref{bj1}, \eqref{bj2},  \eqref{bj3}, and \eqref{AHardy},  that 
  we have  the estimate \eqref{Lem220}, which ends the proof of Lemma \ref{lem22}.
\Box

\subsection{A Lyapunov functional
for equation (\ref{A})} \label{section22}

 With Lemmas \ref{lem21} and \ref{lem22} at hand, we are prepared to establish the first assertion of Proposition \ref{lyap}. Recalling the definition of the functional from \eqref{LL}, namely,
\begin{equation}\label{21juin1bisBnew}
L(w(s),s)=\exp\left(\frac{p+3}{\sqrt{s}}\right) E(w(s),s)+e^{-s},
\end{equation} 
we now provide the proof of estimate \eqref{L01}.
 
 {\it Proof of \eqref{L01}:} 
From  the definition of $E (w(s),s)$ given  in \eqref{EE},   Lemmas  \ref{lem21}, \ref{lem22},
 by using the identity  \eqref{defF123}, 
\begin{eqnarray}\label{23jan12}
\frac1{2s^a}|w|^{p+1}\ln^a(2+ \phi^2 w^2)- e^{-\frac{2(p+1)s}{p-1}}s^{\frac{2a}{p-1}}f(\phi w)=\frac{p-1}{2(p+1)s^a}|w|^{p+1}\ln^a(2+ \phi^2 w^2)\\
-e^{-\frac{2(p+1)s}{p-1}}s^{\frac{2a}{p-1}}\big(f_1(\phi w)+
f_2(\phi w)\big),\no
\end{eqnarray}
we can write for all
  $s'\ge \taa$, 
\begin{align}\label{PP1}
E(w(s'),s') - E(w(s),s) \le &\int_{s}^{s'} \frac{p+3}{2\tau^{\frac32}}E(w(\tau),\tau) {\mathrm{d}}\tau
-\int_{s}^{s'}\Big( \frac3{p-1}-\frac{C_1}{\tau}\Big) \int_{-1}^1 (\partial_{s}w)^2\frac{\rho}{1-y^2}{\mathrm{d}}y{\mathrm{d}}\tau\no\\
&-\int_{s}^{s'}\Big(\frac{p+7}{4}   -\frac{C_1}{\tau^{\frac12}}\Big)\frac1{\tau^{\frac32}}
\int_{-1}^1 (\partial_{s} w)^2\rho {\mathrm{d}}y{\mathrm{d}}\tau\no\\
&-\int_{s}^{s'}\Big( \frac{p-1}{4}   -   \frac{C_1}{\tau^\frac14}\Big)\frac1{\tau^{\frac32}}
\int_{-1}^1(\partial_y w)^2(1-y^2)\rho  {\mathrm{d}}y{\mathrm{d}}\tau \no\\
&-\int_{s}^{s'}\Big(\frac{p+1}{2(p-1)}  -\frac{C_1}{\tau^{\frac14}}\Big)\frac1{\tau^{\frac32}}
\int_{-1}^1w^2 \rho {\mathrm{d}}y{\mathrm{d}}\tau\no\\
&+\int_{s}^{s'}\frac1{\tau^{\frac32}}
\int_{-1}^1w\partial_s w \rho {\mathrm{d}}y{\mathrm{d}}\tau\no\\
&-\int_{s}^{s'}\Big(\frac{p-1}{2(p+1)}-\frac{C_1}{\tau^{\frac14}}\Big)\frac1{\tau^{a+\frac32}}\int_{-1}^1 |w|^{p+1}\ln^a(2+\phi^2 w^2)\rho {\mathrm{d}}y{\mathrm{d}}\tau\no\\
&+\int_{s}^{s'} \chi_1(\tau){\mathrm{d}}\tau +C_0 \int_{s}^{s'} e^{-2\tau}\t,
\end{align}
where $C_1$  stands for some positive constant depending only on $p$ and $a$, and where 
\begin{equation*}
\chi_{1}(\tau)=\frac{ p+3}{2\tau^{\frac32}}
e^{-\frac{2(p+1)\tau}{p-1}}\tau^{\frac{2a}{p-1}}  \ibint  \big(f_1(\phi w) +f_2(\phi w)\big)\rho\y.
\end{equation*}
Applying \eqref{equiv2} and \eqref{equiv3}, we infer that
\begin{equation}\label{xsi8}
\chi_{1}(\tau)\le 
\frac{\tilde C_0}{\tau^{a+\frac52}}\ibint |w|^{p+1}\ln^a(2+\phi^2 w^2)\rho \y+\tilde C_0e^{-2\tau}, \quad \forall \tau \in [s,s'].
\end{equation}
Furthermore, we have 
\begin{equation}\label{09}
\frac1{\tau^{\frac32}}
\int_{-1}^1w\partial_sw \rho {\mathrm{d}}y
\le \frac1{p-1}\int_{-1}^1 (\partial_{s}w)^2\frac{\rho}{1-y^2}\y
+\frac{\tilde C_0}{\tau^3}\ibint w^2\rho \y,\quad \forall \tau \in [s,s'].
\end{equation}
Utilizing \eqref{ss00}, we deduce that for all $\tau \in [s,s']$, the following inequalities hold simultaneously:
\begin{align*}
\frac{1}{p-1} - \frac{C_1}{\tau} &\ge 0, \qquad\qquad
\frac{p+7}{4} - \frac{C_1}{\tau^{1/2}} \ge 0, \\
\frac{p-1}{4} - \frac{C_1}{\tau^{1/4}} &\ge 0, \qquad\qquad
\frac{p+1}{2(p-1)} - \frac{\tilde{C}_0 + C_1}{\tau^{\frac14}} \ge 0, \\
\frac{p-1}{2(p+1)} - \frac{\tilde{C}_0+C_1}{\tau^{1/4}} &\ge 0. \qquad\qquad \qquad\qquad\qquad\qquad
\end{align*}
Consequently, 
for all
  $s'\ge \taa$, we have  
\begin{align}\label{FG}
E(w(s'),s') - E(w(s),s) \le &\int_{s}^{s'} \frac{p+3}{2\tau^{\frac32}}E(w(\tau),\tau) {\mathrm{d}}\tau
-\frac2{p-1}\int_{s}^{s'}  \int_{-1}^1 (\partial_{s}w)^2\frac{\rho}{1-y^2}{\mathrm{d}}y{\mathrm{d}}\tau\no\\
& +C_0 \int_{s}^{s'} e^{-2\tau}\t.
\end{align}
On the other hand, it follows from \eqref{21juin1bisBnew} and a direct computation that
  \begin{equation}\label{239}
  \frac{d}{ds}
L(w(s),s) =\exp\Big(\frac{p+3}{\sqrt{s}}\Big)\Big(\frac{d}{ds}
E(w(s),s)-\frac{p+3}{2s^{\frac32}}E(w(s),s)\Big)-e^{-s}.
\end{equation}

Therefore, combining estimates \eqref{FG} and \eqref{239}, we obtain for all   $s'\ge \taa$,   
\begin{align}\label{BH}
L(w(s'),s') - L(w(s),s) \le &
-\frac2{p-1}\int_{s}^{s'}  \exp\Big(\frac{p+3}{\sqrt{\tau}}\Big) \int_{-1}^1 (\partial_{s}w)^2\frac{\rho}{1-y^2}{\mathrm{d}}y{\mathrm{d}}\tau\no\\
& + \int_{s}^{s'}\Big(C_0\exp\Big(\frac{p+3}{\sqrt{\tau}}\Big)  e^{-2\tau}-e^{-\tau}\Big)\t.
\end{align}
By exploiting \eqref{ss00}, we see that
\begin{equation}\label{theta3}
C_0\exp\left(\frac{p+3}{\sqrt{\tau}}\right)e^{-2\tau} - e^{-\tau} \le 0 \qquad \forall \tau \in [s,s'].
\end{equation}
Using \eqref{theta3} along with the fact that $\exp\left(\frac{p+3}{\sqrt{\tau}}\right)\ge1$, it follows from \eqref{BH} that \eqref{L01} holds. This completes the proof of the first part of Proposition \ref{lyap}. $\Box$

\subsection{Positivity of the Lyapunov functional $L(w(s),s)$}\label{sec23}
The purpose of this subsection is to prove the second part of Proposition \ref{lyap}. 
Specifically, we show that the Lyapunov functional $L(w(s),s)$ satisfies the positivity property given in \eqref{posi}.

{\it Proof of \eqref{posi}}:
 The proof provided in \cite{HZjmaa2020}  addresses the non-characteristic case, where $w = w_{x_0,T_0}$ is derived from a solution $u(x,t)$ of \eqref{gen} at a point $x_0$ satisfying the non-characteristic property. As noted previously, the technical arguments and the Lyapunov functional machinery extend naturally to any case where the cylinder $ (-1,1) \times [-\log T_0, +\infty)$ remains strictly within the interior of the domain of definition of $w$. To conclude the proof for a characteristic point $x_0$, we proceed as follows:

For each $n \in \mathbb{N}$, let us introduce the rescaled function $w_n(z, \sigma)$ defined for all\\
 $\sigma \ge -\ln \left(T_0-\frac{1}{n}\right)$ and $\vert{}z\vert{} < 1 + \frac{e^{\sigma}}{n}$ by:
$${w}_n(z,\sigma)=w_{x_0,T_0-\frac1{n} }(z,\sigma). $$From
(\ref{scaling}), we see that for all $(y,s)\in B\times
[-\ln T_0,+\infty)$
\begin{equation}\label{blowupbefore}
 {w}_n(z,\sigma)=\frac{\phi (-\ln (\frac1{n}+e^{-\sigma}))}{\phi(\sigma)}
w( y,s),
\end{equation}
where $\phi$ defined in \eqref{defphi}, 
and
\begin{equation}\label{sigmaz}
y=\frac{z}{1+\frac{e^{\sigma}}n}, \qquad 
s=\sigma-\ln\left(\frac{e^{\sigma}}{n}+1\right).
\end{equation}

For a fixed $n \in \mathbb{N}$, since $w_{n}(z,\sigma)$ is well defined on the domain $\{\vert{}z\vert{}<2\}$ for large $\sigma$, it follows immediately from the definition that $L(w_{n}(\sigma),\sigma) \to 0$ as $\sigma \to \infty$. 
Consequently, invoking the hypothesis that $(w,\partial_{s}w)(-\log T_0) \in H^{1} \times L^{2}(-1,1)$, we deduce that for any $\sigma \in [-\log T_0,\infty)$ and all sufficiently large $n$,
\begin{equation}\label{nn}
0 \le L(w_{n}(\sigma),\sigma)\le L_{0},
\end{equation}
 where $L_0$ is a positive constant independent of $n$.
By applying the estimates (3.72) established in \cite{HZjmaa2020}, we infer that for all $\sigma \in [1-\log T_0,\infty)$ and $n \in \mathbb{N}$,
\begin{equation*}
\int_{\sigma}^{\sigma+1}\int_{-1}^{1}\left((\partial_{y}w_{n})^{2}(1-z^{2})+
e^{-\frac{2(p+1)\tau}{p-1}}\tau^{\frac{2a}{p-1}}   f(\phi(\tau) w_n)
+w_{n}^{2}+(\partial_{s}w_{n})^{2}\right)\rho (z)dz\t \le CL_{0}+C,
\end{equation*}
where $C$ is a positive constant independent of $n$. Then, by the convergence in the energy space, we obtain for all $\delta>0$ and $s\in[-\log T_0,\infty)$,
\begin{equation*}
\int_{s}^{s+1}\int_{|y|<1-\delta}\left((\partial_{y}w)^{2}(1-y^{2})+
e^{-\frac{2(p+1)\tau}{p-1}}\tau^{\frac{2a}{p-1}}   f(\phi(\tau) w)
+w^{2}+(\partial_{s}w)^{2}\right)\rho(y)\y\t \le CL_{0}+C.
\end{equation*}
Letting $\delta\to0$, we get
\begin{equation*}
\int_{s}^{s+1}\int_{|y|<1}\left((\partial_{y}w)^{2}(1-y^{2})+
e^{-\frac{2(p+1)\tau}{p-1}}\tau^{\frac{2a}{p-1}}   f(\phi(\tau) w)
+w^{2}+(\partial_{s}w)^{2}\right)\rho(y)\y \t\le CL_{0}+C.
\end{equation*}
Furthermore, for all $s\in [-\log T_0,\infty)$ and $|y|<1$ for $n$ large (depending on $s$),
\begin{eqnarray*}
\left((\partial_{y}w_{n})^{2}(1-z^{2})+w_{n}^{2}+e^{-\frac{2(p+1)\sigma}{p-1}}\sigma^{\frac{2a}{p-1}}   f(\phi(\sigma) w_n)
+(\partial_{s}w_{n})^{2}\right)(z,\sigma)\rho(z)\nonumber\\
 \le C\left((\partial_{y}w)^{2}(1-y^{2})+w^{2}+e^{-\frac{2(p+1)s}{p-1}}s^{\frac{2a}{p-1}}   f(\phi(s) w)+(\partial_{s}w)^{2}\right)(y,s)\rho(y),
\end{eqnarray*}
where $(z, \sigma)$ and $(y, s)$ are linked by \eqref{sigmaz}. Therefore, 
 by the Lebesgue theorem,
\begin{equation}
\int_{s}^{s+1}L(w_{n}(\tau),\tau )d\tau \rightarrow \int_{s}^{s+1}L(w(\tau),\tau)d\tau,\qquad \forall s \in [-\log T_0, \infty).
\end{equation}
In view of \eqref{nn} and the monotonicity of $L$, statement \eqref{posi} immediately follows. 
\hfill $\Box$

\bigskip

\noindent {\it Proof of Proposition \ref{lyap}.}\\
The result follows by combining the proof of the first part of Proposition \ref{lyap} presented in Subsection \ref{section22} with that of \eqref{posi} given in Subsection \ref{sec23}.
\hfill $\Box$

\medskip

Combining \eqref{L01} and \eqref{posi}, we conclude that the time-average of the following terms remain bounded:
\begin{cors}\label{cor21}
For all  
 $s \geq  -\ln T_0$
\begin{equation}\label{L1212}
-C\leq E(w(s),s)  \leq  C+CE(w(-\ln T_0),-\ln T_0),
\end{equation}
\begin{equation}\label{cor01Asec4}
\int_{s}^{s+1}\int_{-1}^1(\partial_{s} w)^2\frac{\rho}{1-y^2}{\mathrm{d}}y{\mathrm{d}}\tau\leq
 C +CE(w(-\ln T_0),-\ln T_0).
 \end{equation}
\end{cors}

\section{Proof of Theorem \ref{t1}}\label{section3}
This section is devoted to the proof of Theorem \ref{t1}. The argument is organized into three main steps:
\begin{enumerate}
    \item In Subsection \ref{section31}, by exploiting Corollary \ref{cor21}, we derive the boundedness of both the time-averaged $H^1$ norm and the weighted nonlinear term associated with the solution to \eqref{A}. 
    \item In Subsection \ref{section53}, we establish the uniform boundedness of the $H^1 \times L^2$ norm for the solution to \eqref{A}, thereby concluding the proof of the first part of Theorem \ref{t1}.
    \item In Subsection \ref{sec54}, these estimates are subsequently refined to obtain more precise asymptotic information regarding the behavior of the solution.
\end{enumerate}

 \subsection{Bounds for the time average of the
$H^1$ norm and the nonlinear term of $w$ for solution of equation (\ref{A})}\label{section31}
By exploiting the preceding Corollary, the argument proceeds similarly to the methods developed in \cite{MZajm03,MZimrn05}, with the primary differences lying in the treatment of the nonlinear and perturbation terms. More precisely, we establish the following Proposition:
\begin{prop}\label{prop31}
For all    $s\geq 1-\ln T_0$, we have
\begin{equation}\label{bound1}
\int_{s}^{s+1}\!\ibint \big( w^2(y,\tau )+(\partial_y w(y,\tau ))^2(1-y^2) +\frac{( \partial_{s}w(y,\tau ))^2}{1-y^2}\big)\rho \y \t \leq K_{1},
\end{equation}
\begin{equation}\label{bound2}
\frac{1}{s^a}\int^{s+1}_{s}\! \int_{-1}^1   |w|^{p+1}\ln^a(\phi^2w^2+2) \rho {\mathrm{d}}y{\mathrm{d}}\tau
\leq K_{1} .
\end{equation}
\end{prop}
\begin{proof}\\
-Proof of \eqref{bound2}. 
For $s\ge 1-\ln T_0$, let us work with time integrals between $s_1$ and $s_2$ where $s_1\in [s-1,s]$
and $s_2\in [s+1,s+2]$.
 By integrating the expression \eqref{EE}
 of $E(w(s),s)$ in time between $s_1$ and $s_2$, where $s_2>s_1\ge -\ln T_0$, we obtain:
\begin{align}\label{et}
\int_{s_1}^{s_2}\!E(w(\tau),\tau)\t=&\displaystyle\int_{s_1}^{s_2}\iint\!\!\Big ( \frac{1}{2}(\partial_{s}w)^2
+\frac{p+1}{(p-1)^2}w^2
-e^{-\frac{2(p+1)\tau}{p-1}}\tau^{\frac{2a}{p-1}}   f(\phi w)\Big )\rho {\mathrm{d}}y\t\nonumber\\
&+\frac{1}{2}\displaystyle\int_{s_1}^{s_2}\!\!\iint\!\!  (\partial_y w)^2(1-y^2)\rho {\mathrm{d}}y\t-\frac12\int_{s_1}^{s_2}\!\gamma(\tau)\!\!\displaystyle\iint\!w^2\rho {\mathrm{d}}y\t\no\\
&-\int_{s_1}^{s_2}\!\! \frac1{\tau\sqrt{\tau}}\!\!\displaystyle\iint\!\!w\partial_s w\rho {\mathrm{d}}y\t.
\end{align}
By multiplying the equation (\ref{A}) by $w\rho$ and integrating both in time and in space over $(-1,1)\times [s_1,s_2]$  we obtain the following identity, after some
integration by parts:
\begin{eqnarray}\label{et1}
&&\Big [\iint\!\!\Big (w\partial_{s}w+\frac{5-p}{2(p-1)}w^2\Big ) \rho{\mathrm{d}}y\Big ]_{s_1}^{s_2}=
\int_{s_1}^{s_2}\!\!\iint\!\!(\partial_{s}w)^2\rho{\mathrm{d}}y\t\\
&&-\int_{s_1}^{s_2}\!\!\iint\!\! (\partial_y w)^2(1-y^2)\rho{\mathrm{d}}y\t
-\frac{2p+2}{(p-1)^2}\int_{s_1}^{s_2}\!\!\iint\!\!w^2\rho{\mathrm{d}}y\t\nonumber\\
&&+\int_{s_1}^{s_2}\!\frac1{\tau^a}\iint\!\!|w|^{p+1}\ln^a(2+ \phi^2 w^2)\rho{\mathrm{d}}y\t-\frac{8}{p-1}\!\int_{s_1}^{s_2}\!\iint\!\!w\partial_{s}w
\frac{y^2\rho}{1-y^2}{\mathrm{d}}y\t\nonumber\\
&&+2\!\!\int_{s_1}^{s_2}\!\!\iint\!\!y\partial_y w\partial_{s} w \rho\y\t
+\!\!\int_{s_1}^{s_2}\!\!\iint\!\!\gamma(\tau) w^2\rho\y\t
+\frac{2a}{p-1}\int_{s_1}^{s_2}\frac1{\tau}\iint \partial_s w w\rho \y\t.\nonumber\\
&&+\frac{2a}{p-1}\int_{s_1}^{s_2}\frac1{\tau}\iint y\partial_y w w\rho \y\t.\nonumber
\end{eqnarray}
By combining the identities (\ref{et}),  (\ref{et1}) and exploiting  \eqref{23jan12}, we obtain
\begin{eqnarray}\label{et44}
&&\frac{p-1}{2(p+1)}\int_{s_1}^{s_2}\!\frac1{\tau^a}\iint 
|w|^{p+1}\ln^a(2+ \phi^2 w^2)\rho {\mathrm{d}}y\t=\frac12
\Big [\iint\!\!\big (w\partial_{s}w+\frac{5-p}{2(p-1)}w^2\big ) \rho{\mathrm{d}}y\Big ]_{s_1}^{s_2}
\nonumber\\
&&-\int_{s_1}^{s_2}\!\!\iint\!\!(\partial_{s}w)^2\rho{\mathrm{d}}y\t+\int_{s_1}^{s_2}E(w(\tau),\tau)\t
-\!\!\int_{s_1}^{s_2}\!\!\iint\!\!y \partial_y w\partial_{s} w \rho{\mathrm{d}}y\t+\sum_{i=1}^{6}A_i,
\end{eqnarray}
where
\begin{eqnarray*}
&&A_1=\frac{4}{p-1}\int_{s_1}^{s_2}\!\iint\!\!w\partial_{s}w\,\frac{y^2\rho}{1-y^2}{\mathrm{d}}y\t,\qquad
A_2=-\frac{a}{p-1}\int_{s_1}^{s_2}\!\frac1{\tau}\iint\!\!y\partial_yw\,w\,\rho\,{\mathrm{d}}y\t,\\
&&A_3=-\frac{a}{p-1}\int_{s_1}^{s_2}\!\frac1{\tau}\iint\!\!w\partial_sw\,\rho{\mathrm{d}}y\t,\qquad
A_4=\int_{s_1}^{s_2}\!\frac1{\tau\sqrt{\tau}}\iint\!\!w\partial_sw\,\rho{\mathrm{d}}y\t,\\
&&A_5=\int_{s_1}^{s_2}\!\!\iint e^{-\frac{2(p+1)\tau}{p-1}}\tau^{\frac{2a}{p-1}} f_1(\phi w)\rho {\mathrm{d}}y\t,\qquad
A_6=\int_{s_1}^{s_2}\!\!\iint e^{-\frac{2(p+1)\tau}{p-1}}\tau^{\frac{2a}{p-1}} f_2(\phi w)\rho {\mathrm{d}}y\t.
\end{eqnarray*}
We claim that the proof of \eqref{bound2} follows from the lemma below, which allows us to control 
all terms on the right-hand side of \eqref{et44} by the left-hand side:
\begin{lem}\label{g}
 For all  $s \ge 1-\ln T_0$, 
 for all $\nn >0$,
\begin{equation}\label{control300}
\int_{s_1}^{s_2}\!\!\!\iint\!(\partial_yw)^2(1-y^2) \rho{\mathrm{d}}y\t
\le K_2 +2\int_{s_1}^{s_2}\iint\! e^{-\frac{2(p+1)\tau}{p-1}}\tau^{\frac{2a}{p-1}} f(\phi w)\rho{\mathrm{d}}y\t,\qquad
\end{equation}
\begin{equation}\label{control3}
\int_{s_1}^{s_2}\!\!\!\iint\!|y \partial_y w\partial_s w| \rho{\mathrm{d}}y\t
\le \frac{K_2}{\nn} +K_2\nn\int_{s_1}^{s_2}\iint\! e^{-\frac{2(p+1)\tau}{p-1}}\tau^{\frac{2a}{p-1}} f(\phi w)\rho{\mathrm{d}}y\t,\qquad
\end{equation}
\begin{equation}\label{control1}
\sup_{\tau\in [s_1,s_2]}\iint \!\!w^2(y,\tau)\rho{\mathrm{d}}y\le \frac{K_2}{\nn} +K_2\nn\int_{s_1}^{s_2}\iint\! e^{-\frac{2(p+1)\tau}{p-1}}\tau^{\frac{2a}{p-1}} f(\phi w)\rho{\mathrm{d}}y\t,
\end{equation}
\begin{align}\label{control4}
\iint|w\partial_{s}w|\rho {\mathrm{d}}y\le&
\iint (\partial_{s}w)^2\rho {\mathrm{d}}y+
 \frac{K_2}{\nn} +K_2\nn\int_{s_1}^{s_2}\iint\! e^{-\frac{2(p+1)\tau}{p-1}}\tau^{\frac{2a}{p-1}} f(\phi w)\rho{\mathrm{d}}y\t,
\end{align}
\begin{equation}\label{control5}
\iint \big( (\partial_{s}w(y,s_1))^2+(\partial_{s}w(y,s_2))^2\big)\rho {\mathrm{d}}y
\le K_2,
\end{equation}
\begin{equation}\label{control30}
|A_1|
\le \frac{K_2}{\nn} +K_2\nn
\int_{s_1}^{s_2}\iint\! e^{-\frac{2(p+1)\tau}{p-1}}\tau^{\frac{2a}{p-1}} f(\phi w)\rho{\mathrm{d}}y\t,\qquad
\end{equation}
\begin{equation}\label{A19}
|A_2|+|A_3|+|A_4|
\le 
\frac{K_2}{\nn} +(K_2\nn+\frac{C}{s})\int_{s_1}^{s_2}\iint\! e^{-\frac{2(p+1)\tau}{p-1}}\tau^{\frac{2a}{p-1}} f(\phi w)\rho{\mathrm{d}}y\t,
\end{equation}
\begin{equation}\label{A30}
|A_5|+|A_6|
\le C+ \frac{C}{ s}\int_{s_1}^{s_2}\iint\! e^{-\frac{2(p+1)\tau}{p-1}}\tau^{\frac{2a}{p-1}} f(\phi w)\rho{\mathrm{d}}y\t.
\end{equation}
\begin{equation}\label{26control1}
\int_{s_1}^{s_2}\!\iint\! e^{-\frac{2(p+1)\tau}{p-1}}\tau^{\frac{2a}{p-1}} f(\phi w)\rho{\mathrm{d}}y\t
\le C\int_{s_1}^{s_2}
\frac{1}{\tau^a} \!\!\int_{-1}^1 \!\!  |w|^{p+1}\ln^a(\phi^2w^2+2)\rho {\mathrm{d}}y{\mathrm{d}}\tau
+Ce^{-2s},
\end{equation}
\begin{equation}\label{b26control1}
\int_{s_1}^{s_2}
\frac{1}{\tau^a} \!\int_{-1}^1 \!\!  |w|^{p+1}\ln^a(\phi^2w^2+2)\rho {\mathrm{d}}y{\mathrm{d}}\tau
\le C\int_{s_1}^{s_2}\!\iint\! e^{-\frac{2(p+1)\tau}{p-1}}\tau^{\frac{2a}{p-1}} f(\phi w)\rho{\mathrm{d}}y\t
+Ce^{-2s},
\end{equation}
for some 
$s_1$ and $s_2$, such that  $s_1\in [s-1,s]$
and $s_2\in [s+1,s+2]$.
\end{lem}

\bigskip

Indeed, from (\ref{et44}) and  Lemma \ref{g},
we deduce that
\begin{eqnarray}\label{gg}
\int_{s_1}^{s_2}
\frac{1}{\tau^a} \!\int_{-1}^1 \!\!  |w|^{p+1}\ln^a(\phi^2w^2+2)\rho {\mathrm{d}}y{\mathrm{d}}\tau
&\!\le&\! C_3(K_2\nn +\frac{1}{s})
\int_{s_1}^{s_2}\frac{1}{\tau^a} \!\int_{-1}^1 \!\!  |w|^{p+1}\ln^a(\phi^2w^2+2)\rho {\mathrm{d}}y{\mathrm{d}}\tau\no\\
&&+C_3 \frac{K_2}{\nn} +C_3e^{-2s}.
\end{eqnarray}
Now, we use \eqref{ss00} to write
$\frac{C_3}{s}\le \frac{C_3}{1 - \ln(T(x_0))} \le \frac{1}{4}.$
Furthermore, by choosing $\nu_0$ small enough such that $C_3 K_2 \nu_0 \le \frac{1}{4}$. Therefore, \eqref{gg}
implies
\begin{align*}
\int_{s_1}^{s_2}
\frac{1}{\tau^a} \!\int_{-1}^1 \!\!  |w|^{p+1}\ln^a(\phi^2w^2+2)\rho {\mathrm{d}}y{\mathrm{d}}\tau
\le C{K_2}, \qquad \forall 
s\ge 1-\ln T_0.
\end{align*}
Since $[s,s+1]\subset  [s_1,s_2]$, and $s-s_1\le 1$,
 we derive  
\begin{equation}\label{feb19alphabnvb}
\frac{1}{s^a}\int_{s}^{s+1}\! \int_{-1}^1   |w|^{p+1}\ln^a(\phi^2w^2+2)\rho {\mathrm{d}}y{\mathrm{d}}\tau
\leq K_3.
\end{equation}
 showing that estimate \eqref{bound2} holds.\\

It remains to prove Lemma \ref{g}.\\

  {\it{Proof of  Lemma \ref{g}:}}
 Thanks to \eqref{L1212}, and   \eqref{cor01Asec4}, we can adapt with no difficulty
 the proof  in the unperturbed case \cite{MZajm03,MZimrn05}  (up to some very minor changes),
in order to get   the proof of the estimates \eqref{control300}, \eqref{control3}, \eqref{control1},  \eqref{control4},   \eqref{control5}, and \eqref{A19}. Furthermore,
by using   \eqref{cor01Asec4}, we conclude 
\begin{equation}\label{aa1}
|A_1|\le \frac{ K_{3}}{\nn} +\frac{C\nn}{s^2}\!\int_{s_1}^{s_2}\!\iint\!\!w^2
\frac{y^2\rho}{1-y^2}{\mathrm{d}}y\t.
\end{equation}
 Also, by using   the   Hardy-Sobolev inequality   \eqref{AHardy},  and the estimates  \eqref{control300}, and \eqref{control1},     we  conclude 
\eqref{control30}.\\
Finally,
the results \eqref{A30}, \eqref{26control1}, and \eqref{b26control1},
 follows  immediately from 
 \eqref{equiv1},    \eqref{equiv2}  and  \eqref{equiv3}, 
This concludes the proof of Lemma \ref{g} and the estimate \eqref{bound2}  too.\\

-Proof of \eqref{bound1}.
Let us use \eqref{control300},
and \eqref{control1},       to deduce 
\begin{equation}\label{lp1p}
\int_{s}^{s+1}\!\ibint \big( w^2+(\partial_yw)^2(1-y^2) \big)\rho\y \t \leq K_4 .
\end{equation}
By exploiting \eqref{lp1p}, and 
\eqref{cor01Asec4}, we get \eqref{bound1}.
This concludes the proof of  Proposition \ref{prop31}.

\Box

\end{proof}

\subsection{
Uniform boundedness in $H^1 \times L^2$ for the solution of (\ref{A})}\label{section53}
The boundedness of the $H^1(B)$ norm for the solution to \eqref{A} follows from the boundedness of the functional $E(w(s),s)$ established in \eqref{EE}, together with Proposition \ref{prop31}, 
 combining  the strategy developed in  \cite{MZajm03,MZimrn05, MZjfa07}. The perturbative and nonlinear terms are controlled using the estimates \eqref{A30}, \eqref{26control1}, \eqref{b26control1}, and \eqref{equiv6}. 
More precisely, in this subsection, we establish the proof of the first estimate of Theorem \ref{t1}, namely estimate
 \eqref{30}. To this end, by virtue of Lemma \ref{lemH} and Proposition \ref{prop31}, we first prove the following result:
 \begin{lem} \label{lem33}
 For all $s \ge \frac{3}{2}-\log T_0 $,
\begin{equation}\label{po}
\int_{-1}^{1} \frac{w(y,s)^2}{1 - y^2} \rho(y) dy \le C.
\end{equation}
\end{lem}
\begin{proof}
From the Hardy--Sobolev estimate \eqref{AHardy}
and \eqref{bound1}, we deduce that  for all $ s \ge   \frac{3}{2}-\ln T_0$,
\begin{equation}\label{hardy5}
\int_{s-1/2}^{s+1/2} \int_{-1}^{1} \frac{w(y,\tau)^2}{1 - y^2} \rho(y) dy \t \le C_2.
\end{equation} Thus, there is $s_3(s) \in [s - \frac{1}{2}, s]$ such that 
$ \int_{-1}^{1} \frac{w(y,s_3)^2}{1 - y^2} \rho(y) dy \le 2C_2$. We then have from \eqref{bound1} and \eqref{hardy5}:
\begin{align}\label{65}
\int_{-1}^{1} \frac{w(y,s)^2}{1 - y^2} \rho(y) dy &= \int_{-1}^{1} \frac{w(y,s_3)^2}{1 - y^2} \rho(y) dy + 2 \int_{s_3}^{s} \int_{-1}^{1} \frac{w(y,\tau) \partial_s w(y,\tau)}{1 - y^2} \rho(y) dy \t \nonumber\\
&\le C + \int_{s_3}^{s} \int_{-1}^{1} \frac{w^2(y,\tau)}{1 - y^2} \rho(y) dy \t + \int_{s_3}^{s} \int_{-1}^{1} \frac{(\partial_s w)^2(y,\tau)}{1 - y^2} \rho(y) dy \t.
\end{align}
By \eqref{bound1}, \eqref{hardy5} and \eqref{65}, we deduce \eqref{po}. This concludes the proof of 
Lemma \ref{lem33}.
{\ }\Box
\end{proof}
\medskip

With the above lemma in hand, we proceed in two steps to prove the first estimate of Theorem~\ref{t1}, namely \eqref{30}.

\medskip
\noindent
\textbf{ Step 1:  Control of  $w$ in $L_{\rho}^{p+1}(B)$}

\medskip

We claim the following:
\begin{lem}\label{lem34}  For all 
  $s\ge  2-\log T_0,$  we have 
\begin{equation}\label{lp1}
\int_{s}^{s+1}\!\ibint |w(y,\tau)|^{2p}\rho(y)(1-y^2)\y \t \leq K_5.
\end{equation}
\begin{equation}\label{lp0}
\iint\!|w(y,s )|^{p+1}\rho(y){\mathrm{d}}y
\le K_5.
\end{equation}
\end{lem}

\begin{proof}\\
-Proof of \eqref{lp1}.
For all $s\ge  2-\log T_0$,  we introduce the following change of variables to transform the problem from the bounded interval $(-1, 1)$ to the real line $\mathbb{R}$:
\begin{equation*}
\xi = \frac{1}{2} \log\left(\frac{1+y}{1-y}\right) \quad (\text{equivalently } y = \tanh \xi),
\end{equation*}
and define the rescaled function
\begin{equation*}
\overline{w}(\xi,s) = w(y,s) (1-y^2)^{\frac{1}{p-1}}.
\end{equation*}
This transformation maps the weighted integrals on $(-1, 1)$ into standard integral   forms on $\mathbb{R}$:
\begin{equation}\label{r1}
\int_s^{s+1}\int_{-1}^{1} |w(y,\tau)|^{2p} \rho(y) (1-y^2) \, dy \, d\tau = \int_s^{s+1} \int_{\mathbb{R}} |\overline{w}(\xi,\tau)|^{2p} \, d\xi \, d\tau.
\end{equation}
By applying the Sobolev embedding theorem in two dimensions (space $\xi$ and time $\tau$), we deduce that
\begin{equation}\label{r2}
\int_s^{s+1} \int_{\mathbb{R}} |\overline{w}(\xi,\tau)|^{2p} \, d\xi \, d\tau \leq C \left(\int_s^{s+1} \!\!\int_{\mathbb{R}}\! \left( \overline{w}(\xi,\tau))^2 + (\partial_{\xi}\overline{w}(\xi,\tau))^2 + (\partial_s\overline{w}(\xi,\tau) )^2\right) d\xi \, d\tau \right)^{p}.
\end{equation}
Straightforward calculations, together with the Hardy inequality \eqref{hardy1} and the uniform bound \eqref{bound1}, imply that
\begin{equation}\label{r3}
\int_s^{s+1} \!\int_{\mathbb{R}} \!\overline{w}^2(\xi,\tau) \, d\xi \, d\tau = \int_s^{s+1}\! \int_{-1}^{1}\! w^2(y,\tau) \frac{\rho(y)}{1-y^2} \, dy \, d\tau \le C \int_s^{s+1}\! \|w(\cdot,\tau)\|_{\mathcal{H}_0}^2 \, d\tau \le CK_1.
\end{equation}
Furthermore, a direct computation for the spatial derivative yields
\begin{eqnarray}\label{kl}
\int_s^{s+1}\!\int_{\mathbb{R}} (\partial_{\xi}\overline{w}(\xi,\tau))^2 \, d\xi \, d\tau &\le & C  \int_s^{s+1} \!\int_{-1}^{1} \!(\partial_yw(y,\tau))^2 (1-y^2) \rho(y) \, dy \, d\tau\nonumber
 \\ &&+C \int_s^{s+1}\int_{-1}^{1} \frac{w^2(y,\tau) \rho(y)}{1-y^2} \, dy \, d\tau .
\end{eqnarray}
Thus, by virtue of \eqref{r3} and \eqref{kl}, and taking into account \eqref{bound1}, it follows that
\begin{equation}\label{r4}
\int_s^{s+1}\int_{\mathbb{R}}( \partial_{\xi}\overline{w}(\xi,\tau))^2 \, d\xi \, d\tau \le C \int_s^{s+1} \|w(\cdot,\tau)\|_{\mathcal{H}_0}^2 \, d\tau \le CK_1.
\end{equation}
Finally, considering the time derivative and again exploiting \eqref{bound1}, we obtain
\begin{equation}\label{r5}
\int_s^{s+1}\int_{\mathbb{R}} (\partial_s\overline{w}(\xi,\tau))^2 \, d\xi \, d\tau = \int_s^{s+1}\int_{-1}^{1} \frac{(\partial_sw(y,\tau))^2 \rho(y)}{1-y^2} \, dy \, d\tau \le CK_1.
\end{equation}
By combining the estimates \eqref{r1}--\eqref{r5}, we conclude that \eqref{lp1} holds.

-Proof of (\ref{lp0}).
For all $s\ge  2-\log T_0$, 
  using the mean value theorem, we derive the existence of $\sigma(s)\in [s,s+1]$ such that
\begin{equation}\label{s1}
\iint\!|w(y,\sigma (s))|^{p+1}\rho(y){\mathrm{d}}y
=\int_{s}^{s+1}\iint\!|w(y,\tau
)|^{p+1}\rho(y){\mathrm{d}}y{\mathrm{d}}\tau.
\end{equation}
This transformation maps the weighted integrals on $(-1, 1)$ into standard integral   forms on $\mathbb{R}$:
\begin{equation}\label{r1b}
\int_s^{s+1}\int_{-1}^{1} |w(y,\tau)|^{p+1} \rho(y)  \, dy \, d\tau = \int_s^{s+1} \int_{\mathbb{R}} |\overline{w}(\xi,\tau)|^{p+1} \, d\xi \, d\tau.
\end{equation}
By applying the Sobolev embedding theorem in two dimensions (space $\xi$ and time $\tau$), we deduce that
\begin{equation}\label{r2b}
\int_s^{s+1} \!\int_{\mathbb{R}}\! |\overline{w}(\xi,\tau)|^{p+1} \, d\xi \, d\tau \leq C \left(\int_s^{s+1} \!\!\int_{\mathbb{R}} \!\left( \overline{w}^2(\xi,\tau) + (\partial_{\xi}\overline{w}(\xi,\tau))^2 + (\partial_s\overline{w}(\xi,\tau))^2 \right) d\xi \, d\tau \right)^{\frac{p+1}2}.
\end{equation}
By combining the estimates \eqref{r1b}, \eqref{r2b}, \eqref{r3}, \eqref{r4} and \eqref{r5}, we conclude that 
\begin{equation}\label{s1b}
\int_s^{s+1}\iint\!|w(y,\tau)|^{p+1}\rho(y){\mathrm{d}}y\t
\le CK_1.
\end{equation}
 Let us introduce the following identity 
\begin{align}\label{s2}
\iint\!|w(y,s )|^{p+1}\rho(y){\mathrm{d}}y=&\iint\!|w(y,\sigma (s))|^{p+1}\rho(y){\mathrm{d}}y  +
\int_{\sigma (s)}^{s}\frac{d}{d \tau}\iint |w(y,\tau)|^{p+1}\rho(y){\mathrm{d}}y{\mathrm{d}}\tau.
\end{align}
Combining (\ref{s1}), \eqref{s2} and   the fact that  $xy\le x^2+y^{2},$  we write  for all $s\ge  2-\log T_0,$ 
\begin{align}\label{janv121}
\iint\!|w(y,s )|^{p+1}\rho(y){\mathrm{d}}y
\le& \int_{s}^{s+1}\!\!\iint\!|w(y,\tau
)|^{p+1}\rho(y){\mathrm{d}}y{\mathrm{d}}\tau
 + \int_{s}^{s+1}\!\!\!\iint\!\!|w(y,\tau)|^{2p}(1-y^2)\rho(y){\mathrm{d}}y{\mathrm{d}}\tau\no
 \\
&+C\int_{s}^{s+1}\!\!\iint\!(\partial_s
w(y,\tau))^2\frac{\rho(y)}{1-y^2}{\mathrm{d}}y{\mathrm{d}}\tau.
\end{align}
Thanks to   (\ref{janv121}),   \eqref{s1b},  \eqref{bound1}, and     \eqref{lp1}  we obtain easily 
\eqref{lp0}.
This concludes  the Lemma \ref{lem34}. 
{\ }\Box
\end{proof}

\medskip
\noindent
\textbf{Step 2: Boundedness of individual terms in $E_0(w(s),s)$.}
\medskip
We now claim that each term constituting the Lyapunov functional $E_0(w(s),s)$ is bounded separately.
\begin{lem} \label{lem34b} For all $s \ge 2-\ln T_0$,
      \begin{equation}\label{31}
       \frac1{s^a}  \int_{-1}^{1}  |w(s)|^{p+1}\ln^a (2+\phi^2 w^2)  \rho(y) \, dy
       +\iint\! e^{-\frac{2(p+1)s}{p-1}}s^{\frac{2a}{p-1}} f(\phi w)\rho(y){\mathrm{d}}y
   \le CK_5,
   \end{equation}
   and
      \begin{equation}\label{32}
    \int_{-1}^{1} \left( (\partial_s w(s))^2  +(\partial_y w(s))^2 (1 - y^2) + \frac{w(s)^2}{1-y^2}  \right) \rho(y)  dy \le CK_5. 
    \end{equation}
\end{lem}
\begin{proof}\\
-Proof of \eqref{31}. By using \eqref{lp0},  \eqref{equiv1}, and \eqref{equiv6}, we deduce \eqref{31}.\\
-Proof of \eqref{32}.
By the definition of the functional $L(w(s), s)$ in \eqref{LL}, the estimates established in Proposition \ref{lyap}, and estimates \eqref{31}, we have
\begin{equation}\label{32bis}
\int_{-1}^{1} \left( (\partial_s w(y,s))^2 + (\partial_y w(y,s))^2 (1 - y^2) + w(y,s)^2 \right) \rho(y)  dy \le CK_5.
\end{equation}
By virtue of Lemma \ref{lem33} and the estimate \eqref{32bis}, we obtain \eqref{32}. This concludes  the Lemma \ref{lem34b}. 
{\ }\Box
\end{proof}
Applying Lemma \ref{lem34b} with $T_0=T(x_0)$ completes the proof of the first part of Theorem \ref{t1}.

\subsection{Convergence of positive and negative parts in terms of \texorpdfstring{$E_0(w)$}{E\_0(w)}}
\label{sec54}
    For simplicity, we introduce the notation
\begin{align}
\label{Eplus}
E_{+}(w(s),s) &= \int_{-1}^{1} \left( \frac{1}{2} (\partial_y w)^2 (1 - y^2) + \frac{1}{2} (\partial_s w)^2 + \frac{p+1}{(p-1)^2} w^2 \right) \rho(y) \, dy, \\
\label{Emoins}
E_{-}(w(s),s) &= \frac{1}{s^a} \int_{-1}^{1} |w(y,s)|^{p+1} \ln^a\left(2 + \phi^2(s)\, w^2(y,s)\right) \rho(y) \, dy.
\end{align}

To complete the proof of Theorem \ref{t1}, we first state and prove two preliminary lemmas. The first lemma is given below.
\begin{lem}\label{lem37}
There exists $E_\infty \ge 0$ such that the following assertions hold:
\begin{enumerate}
    \item[(i)] 
    \begin{equation}\label{lim0}
    E_0(w(s),s) \to E_\infty \quad \text{as } s \to \infty.
    \end{equation}
    \item[(ii)] 
    \begin{equation}\label{lim1} 
    \int_s^{+\infty} \int_{-1}^1 \frac{(\partial_s w(y, \tau))^2}{1 - y^2} \rho(y) \, dy \, d\tau \to 0 \quad \text{as } s \to \infty.
    \end{equation}    
    \item[(iii)] \textbf{\textup{(Refined Estimates)}}
    There exist functions $s_4(s)$ and $s_5(s)$ defined for $s \ge 2-\log T_0$ such that
    $  |s_4(s) - s| + |s_5(s) - s| \to 0 \quad \text{as } s \to \infty,$
    and the following energy partition holds as $s \to \infty$:
    \begin{equation}\label{lim2}
        \int_{s_4(s)}^{s_5(s)+1}    E_{-}(w(\tau),\tau)     \, d\tau
 \to \frac{2p+2}{p-1} E_\infty,
    \end{equation}
    \begin{equation}\label{lim3}
    \int_{s_4(s)}^{s_5(s)+1}    E_{+}(w(\tau),\tau)     \, d\tau \to \frac{p+1}{p-1} E_\infty.
    \end{equation}
\end{enumerate}
\end{lem}

\begin{proof}\\

-Proof of \eqref{lim0}. From the properties of the functional $L$ established in Proposition \ref{lyap}, specifically its monotonicity and the fact that it is bounded below by zero, we deduce the existence of a constant $L_{\infty} \ge 0$ such that
\begin{equation}\lim_{s \to \infty} L(w(s), s) = L_{\infty}.
\end{equation}
Therefore, by exploiting \eqref{31} and \eqref{32}, we deduce \eqref{lim0}, where $E
_{\infty}=L_{\infty}\ge0$.

-Proof of \eqref{lim1}: By exploiting \eqref{L01} and the positivity of $L$, we obtain\begin{equation}\label{L0bis1}\int_{s}^{s'} \int_{-1}^1 (\partial_{\tau} w(y, \tau))^2 \frac{\rho}{1-y^2} dy d\tau \le C L(w(s), s), \quad \forall s'\ge s.
\end{equation}
The convergence of the integral on the left-hand side as $s' \to \infty$ directly implies
 \eqref{lim1}.

\noindent{\it Proof of \eqref{lim2} and \eqref{lim3}.} We apply \eqref{et44} between $s_1=s_4(s)$ and $s_2=s_5(s)+1$, which are chosen as follows.
Let $\eta(s) = \int_{s}^{s+1} \int_{-1}^{1} (\partial_s w(y, \tau))^2 \frac{\rho(y)}{1-y^2} \, dy d\tau$. 
Since \eqref{lim1} holds, we have  $\eta(s) \to 0$ as $s \to \infty$. 
Let $s_4(s) \in [s, s + \sqrt{\eta(s)}]$ such that 
\begin{equation}\label{00}
\int_{-1}^1 \partial_s w(s_4(s))^2 \frac{\rho}{1-y^2} \y = \frac{1}{\sqrt{\eta(s)}} \int_s^{s+\sqrt{\eta(s)}} \int_{-1}^1 \partial_s w^2 \frac{\rho}{1-y^2} \y\le \frac{\eta(s)}{\sqrt{\eta(s)}} \to 0.
\end{equation} 
Moreover, we set $s_5(s)=s_4(s+1)-1$, so that $s_5(s)+1=s_4(s+1)$ and, by \eqref{00} at $s+1$,
\begin{equation}\label{001}
 \int_{-1}^{1} (\partial_s w(s_5(s)+1))^2 \frac{\rho}{1-y^2} \to 0 \quad \text{as } s \to \infty.
\end{equation} 
Using \eqref{00}, \eqref{001}, \eqref{32}, and \eqref{lim1} 
and $E(w(s),s) \to E_{\infty}$ and $|s_4(s) - s| + |s_5(s) - s| \to 0$, we conclude that 
$\int_{s_4}^{s_5+1} E(w(s)) \, ds \to E_{\infty}$, providing \eqref{lim2}. Finally, the definition of $E(w(s))$ yields the partition for the remaining terms \eqref{lim3}.
This concludes  the Lemma \ref{lem37}. 
{\ }\Box
\end{proof}
We now state and prove the second preliminary lemma
\begin{lem}
\label{lem38}For all $s \ge 5-\log T_0$,
 $\tau \in [s-3, s+10]$, the function $E_+(w(\tau),\tau)$ satisfies the estimate
\begin{align}
\left|  E_+(w(\tau),\tau) -E_+(w(s),s) \right| &\le C \int^{\tau}_{s} \int_{-1}^{1} \frac{(\partial_s w)^2}{1 - y^2} \rho(y) \, dy \, ds'\label{poi} \\
&\quad + C \int^{\tau}_{s} \left( \int_{-1}^{1} \frac{(\partial_s w)^2}{1 - y^2} \rho(y) \, dy \right)^{1/2} ds' + \frac{C}{s}.\nonumber
\end{align}
\end{lem}
\begin{proof}
Multiplying \eqref{A} by $\partial_{s} w$ and integrating over $(-1,1)$ with weight $\rho(y)$, we obtain for all $\tau \ge  s-3$,
\begin{align}
\frac{d}{d\tau} E_+(w(\tau),\tau) =& -\frac{4}{p-1}\int_{-1}^1 (\partial_{s} w)^2 \frac{\rho(y)}{1-y^2} \, dy + \frac{2a}{(p-1)\tau} \int_{-1}^1 (\partial_s w)^2 \rho(y) \, dy \nonumber\\
& + \frac{1}{\tau^a} \int_{-1}^1 \partial_s w |w|^{p-1}w \ln^a(\phi^2 w^2 + 2) \rho(y) \, dy \nonumber\\
& + \frac{2a}{(p-1)\tau} \int_{-1}^1 \partial_s w \, y \partial_y w \rho(y) \, dy + \gamma(\tau) \int_{-1}^1 w \partial_s w \rho(y) \, dy.
\end{align}
Exploiting \eqref{31} and \eqref{32},
 we obtain, for all $\tau \in [s-3, s+10]$,
\begin{align}\label{tty}
\left| \frac{d}{d\tau} E_+(w(\tau),\tau) \right| &\leq C \left( \int_{-1}^1 \frac{(\partial_s w)^2}{1-y^2} \rho(y) \, dy \right)^{\frac{1}{2}} \left( \frac{1}{\tau^{2a}} \int_{-1}^1 |w|^{2p} \ln^{2a}(\phi^2 w^2 + 2)(1-y^2) \rho(y) \, dy \right)^{\frac{1}{2}} \nonumber\\
&\quad + C \int_{-1}^1 (\partial_s w)^2 \frac{\rho(y)}{1-y^2} \, dy + \frac{C}{\tau}.
\end{align}
Using \eqref{hardy3} and \eqref{32}, we deduce that
\begin{equation}\label{Linfty}
\| w(\tau) (1-y^2)^{\frac{1}{p-1}} \|_{L^\infty(-1,1)} \leq C.
\end{equation}
Furthermore, combining \eqref{Linfty}, \eqref{equiv7}, and \eqref{lp0} yields
\begin{equation}\label{Linfty1}
\frac{1}{\tau^{2a}} \int_{-1}^1 |w|^{2p} \ln^{2a}(\phi^2 w^2 + 2)(1-y^2) \rho(y) \, dy \le C.
\end{equation}
Therefore, combining \eqref{tty}, \eqref{Linfty}, and \eqref{Linfty1}, we arrive at
\begin{align*}
\left| \frac{d}{d\tau} E_+(w(\tau),\tau) \right| \le C \int_{-1}^{1} \frac{(\partial_s w)^2}{1 - y^2} \rho(y) \, dy + C \left( \int_{-1}^{1} \frac{(\partial_s w)^2}{1 - y^2} \rho(y) \, dy \right)^{1/2} + \frac{C}{\tau}.
\end{align*}
Integrating over $[s, \tau]$,  for all $\tau \in [s-3, s+10]$, we deduce \eqref{poi}.
This concludes  the Lemma \ref{lem38}. 
{\ }\Box
\end{proof}


With these preliminaries in hand, we can now establish \eqref{33} and \eqref{34}.
\subsubsection*{Proof of \eqref{33} and \eqref{34}}
\begin{proof}
We have from \eqref{lim3} that for all $\epsilon_0 >0$, there is $s_{\epsilon_0} \ge -\log (T_0) + 5$ such that for all $s \ge s_{\epsilon_0}$, we have:
\begin{equation}
\left| \int_{s_4(s-2)}^{s_5(s-2)+1} 
E_{+}(w(\tau),\tau) \t - \frac{p+1}{p-1} E_{\infty} \right| \le \frac{\epsilon_0}{4}, \label{p0}
\end{equation}
\begin{equation}
\int_{s_4(s-2)-1}^{s} \int_{-1}^{1} \frac{\partial_s w(y, \tau)^2}{1 - y^2} \rho \, dy \, \t \le 
\epsilon_0^4, \label{p1}
\end{equation}
and $|s_4(s) - s| + |s_5(s) - s| \le \delta_0(\epsilon_0)$, where 
$\delta_0(\epsilon_0):=\epsilon_0^4$.\\
From \eqref{p0}, we know that for all $s \ge s_{\epsilon_0}$, there is $s_6 \in [s_4(s - 2), s_5(s - 2) + 1]$ such that,
\[ \left| (1 + s_5 - s_4) 
E_{+}(w(s_6),s_6)
 - \frac{p+1}{p-1} E_{\infty} \right| \le \frac{\epsilon_0}{4}. \]
 Note that $s_6\in[s-2,s]$ and $|s_5-s_4|\le2\delta_0=2\epsilon_0^4$; since $E_+(w(s_6),s_6)\le C$ by \eqref{32}, we deduce
\begin{equation}
\left|E_{+}(w(s_6),s_6)
  - \frac{p+1}{p-1} E_{\infty} \right| \le \frac{\epsilon_0}{2} + C_0 \epsilon_0^4. \label{p3}
\end{equation}
Furthermore, by \eqref{poi} and \eqref{p1}, we have 
\begin{equation}\label{yy}
\left|E_+(w(s),s)- E_+(w(s_6),s_6) \right|\le 
 C\epsilon_0^4+C\epsilon_0^2+\frac{C}{s}.
\end{equation}
By combining  \eqref{yy} and \eqref{p3}  and  by taking $\epsilon_0$ small enough, we deduce \eqref{33}.
Finally, the definition of $E(w(s))$ yields the partition for the remaining terms \eqref{32}.
Thus, \eqref{34} is proved. This concludes the remaining  part of Theorem \ref{t1} 
\Box
\end{proof}

 \appendix

 \section{Some elementary lemmas.}
Let the functions $\phi$, $f$, $f_1$, and $f_2$ be defined as in \eqref{defphi}, \eqref{defF}, \eqref{defF2}, and \eqref{defF123}, respectively.
Clearly, we have 
\begin{lem} \ {\  }Let $r>1$,\\ 
\begin{align}
\int_0^u|v|^{r-1}v\ln^{{a}}(2+v^2 )\v  \sim&  \frac{| u|^{r+1}}{r+1}\ln^{{a}}(u^2+2  ),\quad  \text{ as } \;\; |u| \to \infty,
\label{estF0}\\
f(u)  \sim & \frac{| u|^{p+1}}{p+1}\ln^{{a}}(u^2+2  )  \quad  \text{ as } \;\; |u| \to \infty,\label{estF}\\
f_2(u)\sim & C| u|^{p+1}\ln^{{a-2}}(u^2+2) \quad  \text{ as } \;\; |u| \to \infty.\label{estF3}
\end{align}
\end{lem}
\begin{proof} 
 See Lemma A.1 
 in \cite{HZjmaa2020}.
\end{proof}
\Box

Thanks to \eqref{estF0}, \eqref{estF} and \eqref{estF3}, we can
  state and prove  the following estimates:

\bigskip

%

\begin{lem}\label{lemm:esth}
For all  $s \geq 1$,  and $z\in \R$, we have
\begin{align}
 C^{-1}|\phi z|^{p+1} \ln^a(\phi^2 z^2+2)\le C+f(\phi z)\le C +C|\phi z|^{p+1} \ln^a(2+\phi^2 z^2),\label{equiv1}\\
|f_1(\phi z)|\le C+\frac{C}{s}|\phi z|^{p+1} \ln^a(2+\phi^2 z^2),\quad\quad\label{equiv2}\\
|f_2(\phi z)|\le C+\frac{C}{s^2}|\phi z|^{p+1} \ln^a(2+\phi^2 z^2),\quad\quad\label{equiv3}\\
 \frac{1}{s^a} |z|^{p+1} \ln^a(2 + \phi^2 z^2) 
\le  C| z|^{ p+1}+Ce^{-s}, \quad \textrm{if}\quad  a<0,\quad\label{equiv6}\\
\frac{1}{s^{2a}} |z|^{p+1} \ln^{2a}(2 + \phi^2 z^2) 
\le  C| z|^{ p+1}+Ce^{-s}, \quad \textrm{if}\quad  a<0.\quad\label{equiv7}
\end{align}
Furthermore,    there exists $\tilde{S}_0 > 1$ such that for all $s \ge \tilde{S}_0$,  $\varepsilon > 0$, and $z \in \mathbb{R}$, the following estimate holds:
\begin{equation}\label{equiv5}
z^2 \le \frac{\varepsilon}{s^a} |z|^{p+1} \ln^a(2 + \phi^2(s) z^2) + C(\varepsilon), \quad \textrm{if}\quad  a<0.
\end{equation}
\end{lem}

\begin{proof} 
The estimate \eqref{equiv1} follows immediately from \eqref{estF}. The derivation of \eqref{equiv2} and \eqref{equiv3} is obtained by partitioning the domain into $\{z^2\phi \geq 4\}$ and $\{z^2\phi \leq 4\}$ and applying \eqref{estF0}--\eqref{estF3}.
Furthermore, the assumption $a < 0$ yield \eqref{equiv6}  and \eqref{equiv7} via a similar argument. We prove the inequality \eqref{equiv5} by distinguishing two cases regarding the size of $|z|$:
\begin{itemize}
\item \textbf{Case 1:} $|z| \le \phi(s)$: Since the function $\d \xi \mapsto \ln^a(2+\xi)$ is strictly decreasing on $[2,\infty)$, it follows that
\begin{equation}\label{equiv8}
\frac1{s^a} \ln^a(2+\phi^2(s)z^2) \ge \frac1{s^a}\ln^a(2+\phi^4(s)).
\end{equation}
By the asymptotic property $\d \lim_{s \to \infty} \frac{1}{s^a} \ln^a(2 + \phi^4(s)) = \frac{8^a}{(p-1)^a}  > 0$, there exists $\tilde{S}_0 > 1$ such that
\begin{equation}\label{equiv8b}
\frac{1}{s^a} \ln^a(2 + \phi^4(s)) \ge \frac{8^a}{2(p-1)^a}, \quad \forall s \ge \tilde{S}_0.
\end{equation}
Combining \eqref{equiv8} and \eqref{equiv8b}, we obtain
\begin{equation}\label{equiv8c}\frac1{s^a}\ln^a(2+\phi^2(s)z^2) \ge  \frac{8^a}{2(p-1)^a}, \quad \forall s \ge \tilde{S}_0.
\end{equation}
Applying Young’s inequality in the form $z^2 \le \nu |z|^{p+1} + C(\nu)$ with the choice $\nu = \frac{\varepsilon 8^a}{2(p-1)^a}$, we arrive at \eqref{equiv5}.\item \textbf{Case 2:} $|z| \ge \phi(s)$.Using the monotonicity of the logarithm, we have$$\frac1{s^a}\ln^a(2+\phi^2(s)z^2) \ge \frac1{s^a}\ln^a(2+z^4).$$Standard asymptotic analysis implies that for any $\varepsilon > 0$, there exists a constant $C(\varepsilon) > 0$ such that
\begin{equation}\label{aux_est}z^2 \le \varepsilon |z|^{p+1} \ln^a(2 + z^4) + C(\varepsilon).
\end{equation}
Combining \eqref{aux_est} with the observation that $s^{-a} \ge 1$ for $s \ge \tilde{S}_0$, we conclude the proof of \eqref{equiv5}.
\end{itemize}
This concludes the proof of Lemma \ref{lemm:esth}. \Box
\end{proof}

Finally, we
recall the following lemma:
\begin{lem}[A Hardy-Sobolev type inequality]\label{lemH}
There exists a constant $C>0$ such that, for all $h\in\mathcal H_0$, the following inequalities hold:
\begin{align}
\left( \int_{-1}^{1} |h(y)|^2 \frac{\rho(y)}{1-y^2} dy \right)^{1/2} &\leq C \|h\|_{\mathcal{H}_0}, \label{hardy1} \\
\int_{-1}^{1}|h(y)|^{p+1}\rho (y)dy &\leq C \|h\|^{p+1}_{\mathcal{H}_0}, \label{hardy2} \\
\|h(1-y^2)^{\frac{1}{p-1}}\|_{L^\infty(-1,1)} &\leq C \|h\|_{\mathcal{H}_0}. \label{hardy3}
\end{align}
\end{lem}
\begin{proof} 
 See Lemma 2.2 
 in \cite{ MZjfa07}.
 \Box
\end{proof}



\noindent{\bf Address}:\\
 Department of Basic Sciences, Deanship of Preparatory and Supporting Studies,
  Imam Abdulrahman Bin Faisal University
P.O. Box 1982 Dammam, Saudi Arabia.\\
\vspace{-7mm}
\begin{verbatim}
e-mail:  mahamza@iau.edu.sa
\end{verbatim}
Universit\'e Sorbone Paris Nord, LAGA, CNRS (UMR 7539), F-93430 Villetaneuse, France.\\
\vspace{-7mm}
\begin{verbatim}
e-mail: Hatem.Zaag@univ-paris13.fr
\end{verbatim}



\begin{thebibliography}{10}

\bibitem{Apndeta95}
S. Alinhac,  Blow-up for nonlinear hyperbolic equations, in: Progress in Nonlinear Differential Equations and their Applications, vol. 17, Birkh\"auser Boston Inc., Boston, MA, 1995.






\bibitem{CZcpam13}
R. C\^ote, H. Zaag, Construction of a multisoliton blow-up solution to the semilinear wave equation in one space dimension, Comm. Pure Appl. Math. \textbf{66} (2013) 1541--1581.










\bibitem{Hamza2026}
M.A. Hamza, The blow-up rate for a log non-scaling invariant semilinear wave equation in the conformal regime, arXiv preprint arXiv:2604.05202 (2026).


\bibitem{H1}
M.A. Hamza, The blow-up rate for strongly perturbed semilinear wave equations in the conformal regime without a radial assumption, Asymptot. Anal. \textbf{97} (3-4) (2016) 351--378.

\bibitem{omar1}
M.A. Hamza, O. Saidi, The blow-up rate for strongly perturbed semilinear wave equations, J. Dynam. Differential Equations \textbf{26} (2014) 1115--1131.

\bibitem{omar2}
M.A. Hamza, O. Saidi, The blow-up rate for strongly perturbed semilinear wave equations in the conformal case, Math. Phys. Anal. Geom. \textbf{18} (1) (2015) Art. 15.

\bibitem{HZjhde12}
M.A. Hamza, H. Zaag, A Lyapunov functional and blow-up results for a class of perturbations for semilinear wave equations in the critical case, J. Hyperbolic Differ. Equ. \textbf{9} (2012) 195--221.

\bibitem{HZnonl12}
M.A. Hamza, H. Zaag, A Lyapunov functional and blow-up results for a class of perturbed semilinear wave equations, Nonlinearity \textbf{25} (2012) 2759--2773.

\bibitem{HZkg12}
M.A. Hamza, H. Zaag, Blow-up behavior for the Klein-Gordon and other perturbed semilinear wave equations, Bull. Sci. Math. \textbf{137} (8)(2013) 1087--1109.

\bibitem{HZjmaa2020}
M.A. Hamza, H. Zaag, The blow-up rate for a non-scaling invariant semilinear wave equations, J. Math. Anal. Appl. \textbf{483} (2) (2020) 123617.

\bibitem{HZ2022}
M.A. Hamza, H. Zaag, The blow-up rate for a non-scaling invariant semilinear wave equations in higher dimensions, Nonlinear Anal. \textbf{212} (2022) 112440.






\bibitem{MZajm03}
F. Merle, H. Zaag, Determination of the blow-up rate for the semilinear wave equation, Amer. J. Math. \textbf{125} (2003) 1147--1164.

\bibitem{MZimrn05}
F. Merle, H. Zaag, Blow-up rate near the blow-up surface for semilinear wave equations, Internat. Math. Res. Notices \textbf{19} (2005) 1127--1156.

\bibitem{MZma05}
F. Merle, H. Zaag, Determination of the blow-up rate for a critical semilinear wave equation, Math. Ann. \textbf{331} (2005) 395--416.

\bibitem{MZjfa07}
F. Merle, H. Zaag, Existence and universality of the blow-up profile for the semilinear wave equation in one space dimension, J. Funct. Anal. \textbf{253} (2007) 43--121.

\bibitem{MZcmp08}
F. Merle, H. Zaag, Openness of the set of non-characteristic points and regularity of the blow-up curve for the 1D semilinear wave equation, Comm. Math. Phys. \textbf{282} (2008) 55--86.

\bibitem{MZbsm11}
F. Merle, H. Zaag, Blow-up behavior outside the origin for a semilinear wave equation in the radial case, Bull. Sci. Math. \textbf{135} (2011) 353--373.

\bibitem{MZajm11}
F. Merle, H. Zaag, Existence and classification of characteristic points at blow-up for a semilinear wave equation in one space dimension, Amer. J. Math. \textbf{134} (2012) 581--648.

\bibitem{MZisol10}
F. Merle, H. Zaag, Isolatedness of characteristic points for a semilinear wave equation in one space dimension, Duke Math. J. \textbf{161} (2012) 2837--2908.

\bibitem{MZods}
F. Merle, H. Zaag, On the stability of the notion of non-characteristic point and blow-up profile for semilinear wave equations, Comm. Math. Phys. \textbf{333} (3) (2015) 1529--1562.
  
\bibitem{MZods14}
F. Merle, H. Zaag, Dynamics near explicit stationary solutions in similarity variables for solutions of a semilinear wave equation in higher dimensions, Trans. Amer. Math. Soc. \textbf{368} (1) (2016) 27--87.








\bibitem{RZ}
T. Roy, H. Zaag, The blow-up rate for a loglog non-scaling invariant semilinear wave equation, Nonlinear Differ. Equ. Appl. NoDEA \textbf{32} (2025) 115.















\end{thebibliography}
\end{document}